\documentclass{amsart}

\newtheorem{theorem}{Theorem}[section]
\newtheorem{lemma}[theorem]{Lemma}

\newtheorem{corollary}[theorem]{Corollary}

\usepackage{graphicx}

\theoremstyle{definition}

\theoremstyle{remark}

\numberwithin{equation}{section}

\begin{document}

\title[Hyperbolic knots realizing the maximal distance]{On hyperbolic knots realizing the maximal distance between toroidal surgeries}

%    Information for first author
\author{Masakazu Teragaito}
%    Address of record for the research reported here
\address{Department of Mathematics and Mathematics Education, Hiroshima University,
1-1-1 Kagamiyama, Higashi-hiroshima, Japan 739-8524}
%    Current address
%\curraddr{Department of Mathematics and Statistics,
%Case Western Reserve University, Cleveland, Ohio 43403}
\email{teragai@hiroshima-u.ac.jp}
%    \thanks will become a 1st page footnote.
\thanks{
Partially supported by Japan Society for the Promotion of Science,
Grant-in-Aid for Scientific Research (C), 16540071.
}%

%    General info

%\date{January 1, 1994 and, in revised form, June 22, 1994.}

%\dedicatory{This paper is dedicated to our authors.}


\begin{abstract}
For a hyperbolic knot in the $3$-sphere, the distance between toroidal surgeries is at most $5$, except the figure eight knot.
In this paper, we determine all hyperbolic knots that admit two toroidal surgeries with distance $5$.
\end{abstract}

\maketitle

%%%%%%%%%%%%%%%%%%%%%%%%%%%%%%%%%%%%%%%%%%%%%%%%%%%%%%%%%%%%%%%%%%%%%%%%
\section{Introduction}

Let $K$ be a  hyperbolic knot in the $3$-sphere $S^3$ with exterior $E(K)=S^3-\mathrm{Int}\,N(K)$.
For a slope $\gamma$ on $\partial E(K)$, $K(\gamma)$ denotes the manifold obtained by $\gamma$-surgery.
That is, $K(\gamma)=E(K)\cup V$, where $V$ is a solid torus attached to $E(K)$ along $\partial E(K)$ in such a way that
$\gamma$ bounds a meridian disk in $V$.
By Thurston's hyperbolic Dehn surgery theorem \cite{Th}, $K(\gamma)$ is a hyperbolic $3$-manifold, except at most finitely many slopes $r$.
If $K(\gamma)$ fails to be hyperbolic, then it is conjectured that $K(\gamma)$ is a Seifert fibered manifold or a toroidal manifold \cite{G2}.
A $3$-manifold is said to be \textit{toroidal\/} if it contains an essential torus.
If $K(\gamma)$ is toroidal, then the slope $\gamma$ (or the surgery) is said to be \textit{toroidal\/}.
By \cite{GL}, a toroidal slope corresponds to an integer or a half-integer under the standard parameterization of slopes by the set $\mathbb{Q}\cup \{1/0\}$
(see \cite{R}).
In other words, a toroidal slope runs at most twice along the knot.
Eudave-Mu\~{n}oz \cite{EM} constructed an infinite family of hyperbolic knots, denoted by $k(\ell,m,n,p)$ in his notation,
each of which admits a half-integral toroidal surgery.
In fact, Gordon and Luecke \cite{GL3} proved that these are the only knots that admit half-integral toroidal surgery.

Let $\Delta(\alpha,\beta)$ denote the distance between slopes $\alpha$ and $\beta$.
That is, it is the minimal geometric intersection number between $\alpha$ and $\beta$.
As it is well known \cite{G2},
the figure eight knot admits exactly three toroidal slopes $-4$, $0$ and $4$.
Notice that $\Delta(-4,4)=8$.
For convenience, assume that $K$ is not the figure eight knot.
Then the distance between toroidal slopes of $K$ is at most $5$ by \cite{G}, and
one between integral toroidal slopes is at most $4$ by \cite{T}.
Thus we see that if $K$ admits two toroidal slopes at distance $5$ then one slope is half-integral, and therefore
$K$ is one of Eudave-Mu\~{n}oz knots $k(\ell,m,n,p)$.
In fact, any $k(2,-1,n,0)$, for $n\ne 1$, admits two toroidal slopes $25n-16$ and $25n-37/2$ of distance $5$ \cite{EM2}.
We remark that $k(2,-1,n,0)$ is obtained from the trefoil component of the Whitehead sister link (the $(-2,3,8)$-pretzel link)
by doing $-1/(n-1)$-surgery along the unknotted component,
and that $k(2,-1,1,0)$ is the trefoil and $k(2,-1,0,0)$ is the mirror image of the $(-2,3,7)$-pretzel knot.
The purpose of this paper is to prove that these are the only hyperbolic knots that realize distance $5$ between toroidal slopes.

\begin{theorem}\label{thm:main}
Let $K$ be a hyperbolic knot in $S^3$.
If $K$ admits two toroidal slopes $\alpha$ and $\beta$ with $\Delta(\alpha,\beta)=5$, then $K$ is
the Eudave-Mu\~{n}oz knot $k(2,-1,n,0)$ for some integer $n\ne 1$, and $\{\alpha,\beta\}=\{25n-16, 25n-37/2\}$.
\end{theorem}

\begin{figure}[tb]
%\blankbox{1.0\columnwidth}{2in}
\includegraphics*[scale=0.7]{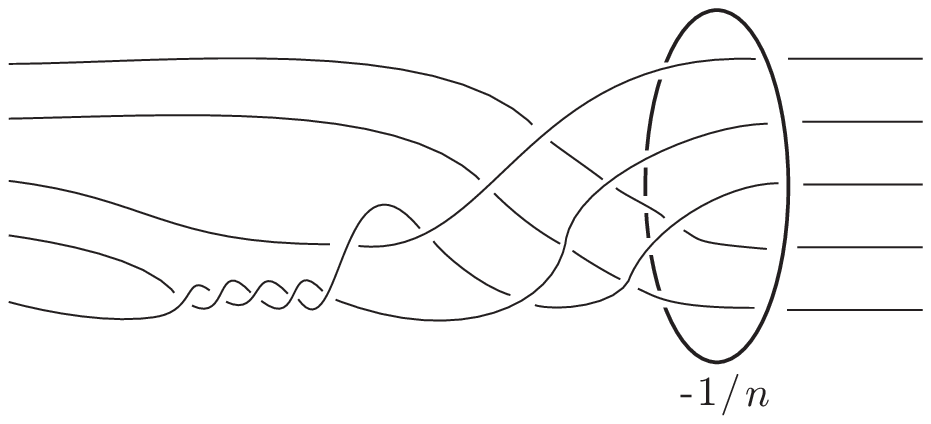}
\caption{}\label{fig:knot}
\end{figure}

In Figure \ref{fig:knot}, the knot $k(2,-1,n,0)$, expressed in a braid form, is obtained after $-1/n$-surgery along the unknotted circle.
This is directly obtained from \cite[Figure 24]{EM}. 
Notice that $k(2,-1,n_1,0)$ and $k(2,-1,n_2,0)$ are not equivalent if $n_1\ne n_2$, since
any knot admits at most one half-integral toroidal slope \cite{GWZ}.

Among the toroidal slopes listed in \cite{EM}, the distance $5$ is realized only by the knots in our Theorem \ref{thm:main}.
If the list of \cite{EM} is complete, then our Theorem would be its consequence.
But it seems to be hard to confirm the completeness, and hence we will take another approach.
As in \cite{G,T}, we consider the labelled intersection graphs on two essential tori in the surgered manifolds.

It is conjectured that a hyperbolic knot in $S^3$ has at most $3$ toroidal slopes \cite{EM} (see also \cite[Problem 1.77 A(5)]{K}).
In \cite{T}, we gave an upper bound $5$ for the number of toroidal slopes.
Theorem \ref{thm:main} implies that the conjecture holds for Eudave-Mu\~{n}oz knots.

\begin{corollary}\label{cor:main2}
Any Eudave-Mu\~{n}oz knot admits at most $3$ toroidal slopes.
\end{corollary}

\begin{proof}
Let $K=k(\ell,m,n,p)$ be an Eudave-Mu\~{n}oz knot.
Then $K$ has the unique non-integral toroidal slope $r$ by \cite{GWZ}.
Let $s=r+1/2$.
Then the slopes $s-1$ and $s$ yield atoroidal Seifert fibered manifolds \cite{EM,EM2}.

If $K$ is not equivalent to $k(2,-1,n,0)$, then the distance between two toroidal slopes is at most $4$ by Theorem \ref{thm:main}.
Hence $\{s-2,r,s+1\}$ gives the set of possible toroidal slopes.
If $K=k(2,-1,n,0)$, then $s=25n-18$ and $s+1$ yields an atoroidal Seifert fibered manifold \cite{EM2}.
Hence $\{s-2,r,s+2\}=\{25n-20,25n-37/2,25n-16\}$ is the set of possible toroidal slopes.
\end{proof}

In Section \ref{sec:prelim}, we collect basic facts for the analysis of the pair of graphs.
In Sections \ref{sec:s4} and \ref{sec:s2}, we will show that there is only one possible pair for the graphs realizing the distance $5$.
Finally, we prove that the pair determines the knot by using the construction of \cite{GL3} in Section \ref{sec:knot}.

%%%%%%%%%%%%%%%%%%%%%%%%%%%%%%%%%%%%%%%%%%%%%%%%%%
\section{Preliminaries}\label{sec:prelim}

Let $K$ be a hyperbolic knot.
Assume that $K$ admits two toroidal slopes $\alpha$ and $\beta$ with distance $5$.
Then one of the slopes is half-integral by \cite{T}.
Hence we assume that $\beta$ is half-integral.
By \cite{GL3}, we know that $K$ is one of Eudave-Mu\~{n}oz knots.
Notice that $K(\alpha)$ and $K(\beta)$ are irreducible by \cite{O,W}.
For the attached solid torus $V_\gamma$ of $K(\gamma)$ where $\gamma\in \{\alpha,\beta\}$, $K_\gamma$ denotes the core of $V_\gamma$.

Let $\widehat{S}$ be an essential torus in $K(\alpha)$.

\begin{lemma}\label{lem:sep}
$\widehat{S}$ is separating in $K(\alpha)$.
\end{lemma}

\begin{proof}
If not, $\alpha=0$ by homological reason.  Thus $K(0)$ contains a non-separating torus, and then $K$ has genus one by \cite{Ga}.
But a genus one knot does not have half-integral toroidal surgery \cite{T2}.
\end{proof}

We can assume that $\widehat{S}$ intersects the attached solid torus $V_\alpha$ of $K(\alpha)$ in a disjoint union of $s$ meridian disks,
$u_1,u_2,\dots,u_s$ numbered successively along $V_\alpha$,
and that $s$ is minimal among all essential tori in $K(\alpha)$.
By Lemma \ref{lem:sep}, $s$ is even.  Let $S=\widehat{S}\cap E(K)$.
Then $S$ is an incompressible and boundary-incompressible,
punctured torus properly embedded in $E(K)$, each of whose boundary components has slope $\alpha$ on $\partial E(K)$.

\begin{lemma}\label{lem:noklein}
$K(\beta)$ does not contain a Klein bottle.
\end{lemma}

\begin{proof}
Since $\beta$ is not integral, this follows from \cite[Theorem 1.3]{GL}.
\end{proof}

Let $\widehat{T}$ be an essential torus in $K(\beta)$.
Then $\widehat{T}$ is separating and we can assume that $\widehat{T}$ intersects the attached solid torus $V_\beta$ in just two meridian disks
$v_1$ and $v_2$ \cite{GL}.  Let $T=\widehat{T}\cap E(K)$.
Then $T$ is an incompressible and boundary-incompressible, twice-punctured torus properly embedded in $E(K)$, where
each component of $\partial T$ has slope $\beta$.

We can assume that $S$ and $T$ intersect transversely, and that no arc component of $S\cap T$ is parallel to
$\partial S$ in $S$ or $\partial T$ in $T$.
Also, we can assume that
no circle component of $S\cap T$ bounds a disk in $S$ or $T$ by the incompressibility of $S$ and $T$, and
that $\partial u_i$ meets $\partial v_j$ in $5$ points for any pair of $i$ and $j$.

If we choose a basis $\{\mu,\alpha\}$ for $H_1(\partial E(K))$, where $\mu$ is the meridian of $K$,
then $\beta$ can be represented by $2\alpha\pm 5 \mu$, since $\Delta(\mu,\beta)=2$ and $\Delta(\alpha,\beta)=5$.
(This fact is called that the jumping number between $\alpha$ and $\beta$ is two in the literature \cite{G,GW}.)

\begin{lemma}\label{lem:jumping}
Let $a_1,a_2,a_3,a_4,a_5$ be the points in $\partial u_i\cap \partial v_j$, numbered
successively along $\partial u_i$.
Then these points appear in the order of $a_1,a_3,a_5,a_2,a_4$ in some direction on $\partial v_j$.
In other words, two points in $\partial u_i\cap \partial v_j$ are adjacent on $\partial u_i$ if and only if
they are not adjacent on $\partial v_j$, and vice-versa.
\end{lemma}

\begin{proof}
This follows from \cite[Lemma 2.10]{GW}.
\end{proof}

In the usual way \cite{CGLS,G}, we obtain two graphs $G_S$ on $\widehat{S}$ and $G_T$ on $\widehat{T}$.
More precisely, $G_S$ has $s$ vertices $u_1,u_2,\dots,u_s$, and $G_T$ has two vertices $v_1$ and $v_2$.
The edges are the arc components of $S\cap T$.
Note that neither graph has a trivial loop.
An endpoint of an edge $e$ in $G_S$ at $u_i$ has label $j$ if it is in $\partial u_i\cap \partial v_j$.
Thus the pair of labels $1$ and $2$ appears $5$ times around $u_i$.
Similarly, the endpoints of edges in $G_T$ are labelled, and then the sequence of labels $1,2,\dots,s$ appears $5$ times around $v_j$.
An edge is called a \textit{$\{i,j\}$-edge\/} if it has labels $i$ and $j$ at its endpoints.

Let $G=G_S$ or $G_T$.
An edge of $G$ is said to be \textit{positive\/} if it connects
two vertices (possibly, the same vertex) with the same parity.  Otherwise, it is said to be \textit{negative}.
We denote by $G^+$ the subgraph of $G$ consisting of all vertices and positive edges of $G$.

A subgraph $H$ of $G$ on a torus $F$ is said to have an \textit{annulus support\/} if there is an annulus $A$ on $F$
such that $H\subset \mathrm{Int}\,A$ and a core of $A$ is essential on $F$, and there is no disk on $F$ containing $H$.

A cycle in $G$ is called a \textit{Scharlemann cycle\/} if it bounds a disk face of $G$ and all the edges are positive
$\{i,i+1\}$-edges for some $i$.
The pair of labels $\{i,i+1\}$ is called the \textit{label pair\/} of the Scharlemann cycle.
In particular, a Scharlemann cycle consisting of two edges is called an \textit{$S$-cycle} for short.

%%%%%%%%%%%%
\begin{lemma}\label{lem:common}
\begin{itemize}
\item[(1)] No two edges can be parallel in both graphs.
\item[(2)] \textup{(The parity rule)} An edge is positive in one graph if and only if it is negative in the other.
\item[(3)] If $\rho$ is a Scharlemann cycle in $G_S$ \textup{(}resp. $G_T$\textup{)},
then its edges cannot lie in a disk on $\widehat{T}$ \textup{(}resp. $\widehat{S}$\textup{)}.
In particular, if $\rho$ is an $S$-cycle, then its edges form a cycle with an annulus support in the other graph.
\end{itemize}
\end{lemma}

\begin{proof}
(1) is \cite[Lemma 2.1]{G}.
(2) can be found in \cite[p.279]{CGLS}.
(3) is \cite[Lemma 3.1]{GL} (Recall that $K(\alpha)$ and $K(\beta)$ are irreducible.)
\end{proof}

\begin{figure}[tb]
%\blankbox{1.0\columnwidth}{2in}
\includegraphics*[scale=0.3]{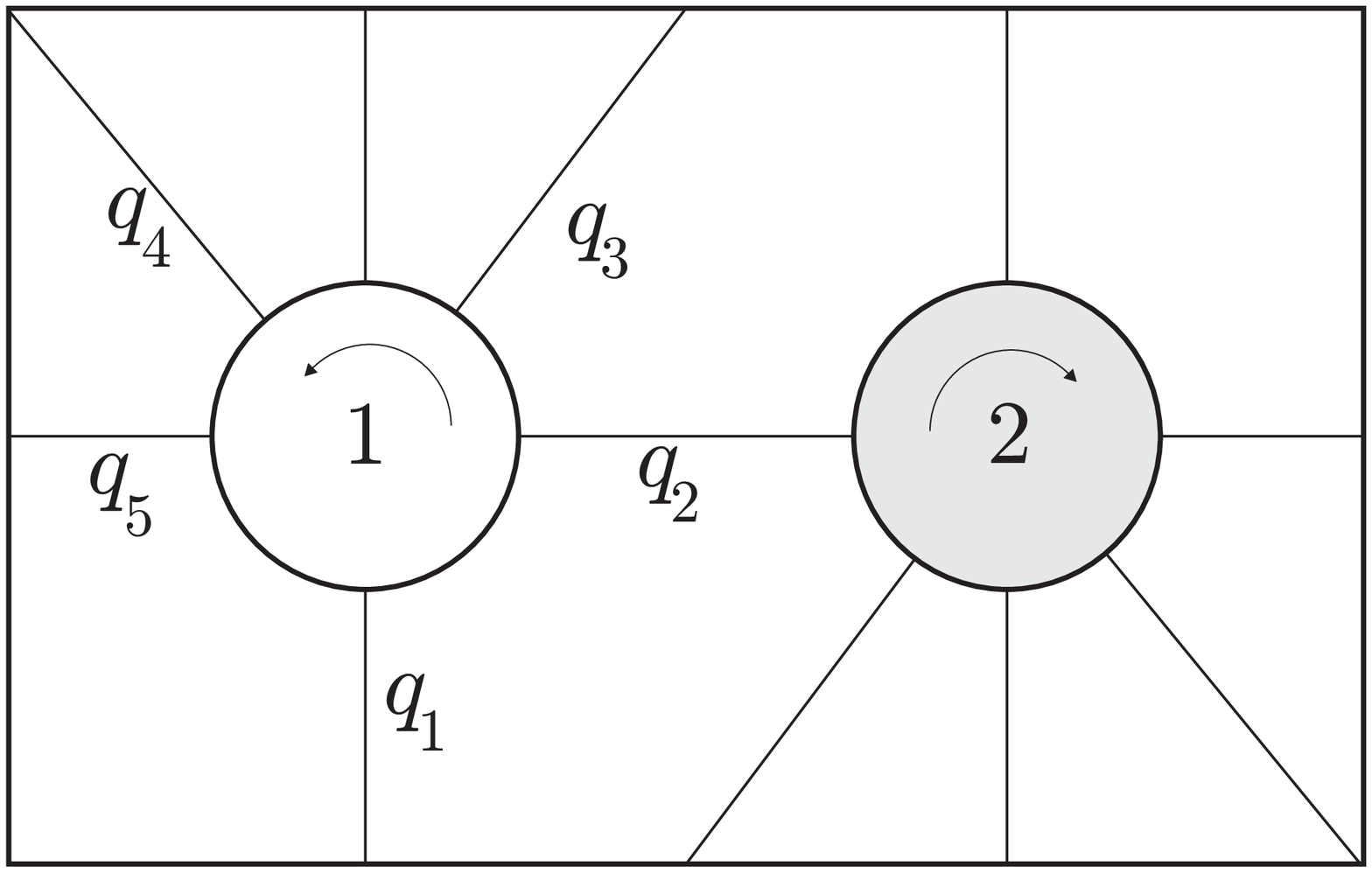}
\caption{}\label{fig:t2}
\end{figure}

%%%%%%%%%%%%%%%%%%%%%%%%%%%%%
The \textit{reduced graph\/} $\overline{G}$ of $G$ is obtained from $G$ by amalgamating each family of parallel edges into a single edge.
For an edge $e$ of $\overline{G}$, the \textit{weight\/} of $e$ is the number of edges in the corresponding family of parallel edges in $G$.
In particular, $\overline{G}_T$ is a subgraph of the graph in Figure \ref{fig:t2}, where
the sides of the rectangle are identified to form the torus $\widehat{T}$ in the usual way (\cite[Lemma 5.2]{G}).
Also, $q_i$ denotes the weight of an edge.
Notice that $v_1$ and $v_2$ are incident to the same number of loops in $G_T$.
Since $G_T$ is determined by non-negative integers $q_i$, we say $G_T\cong G(q_1,q_2,q_3,q_4,q_5)$.
Let $Q_i$ denote the family of parallel edges in $G_T$ with weight $q_i$ for $i=2,3,4,5$.

\begin{lemma}\label{lem:GS}
$G_S$ satisfies the following.
\begin{itemize}
\item[(1)] Any family of parallel positive edges contains at most $3$ edges.
\item[(2)] Any family of parallel negative edges contains at most two edges.
\end{itemize}
\end{lemma}

\begin{proof}
(1) If there are $4$ parallel positive edges in $G_S$, then
there are two bigons among them which lie on the same side of $\widehat{T}$.
By Lemma \ref{lem:common}(1), these $4$ edges belong to mutually distinct families of parallel negative edges in $G_T$.
But this implies that $K(\beta)$ contains a Klein bottle (see the proof of \cite[Lemma 5.2]{G-L2}), which contradicts Lemma \ref{lem:noklein}.

(2) Any negative edge in $G_S$ corresponds to a loop at $v_1$ or $v_2$ in $G_T$ by the parity rule.
The result follows from that any two loops at $v_i$ are parallel and Lemma \ref{lem:common}(1).
\end{proof}

\begin{lemma}\label{lem:GT}
$G_T$ satisfies the following.
\begin{itemize}
\item[(1)] If $s=4$, then there are no two $S$-cycles with disjoint label pairs.
\item[(2)] If $s\ge 4$, then any family of parallel positive edges contains at most $s/2+1$ edges.
Moreover, if it contains $s/2+1$ edges, then two adjacent edges on one end form an $S$-cycle.
\item[(3)] Any family of parallel negative edges contains at most $s$ edges.
\item[(4)] The faces of two $S$-cycles with disjoint label pairs lie on the same side of $\widehat{S}$.
\end{itemize}
\end{lemma}

\begin{proof}
(1) Let $\rho_1$ and $\rho_2$ be $S$-cycles with disk faces $f_1$ and $f_2$, and label pairs $\{1,2\}$ and $\{3,4\}$, say.
Let $H_{12}$ and $H_{34}$ be the parts of $V_\alpha$ between $u_1$ and $u_2$, $u_3$ and $u_4$, respectively.
Then shrinking $V_{12}$ radially to its core in $V_{12}\cup f_1$ gives a M\"{o}bius band $B_1$ such that $\partial B_1$ is the loop on $\widehat{S}$
formed by the edges of $\rho_1$.
Similarly, we obtain another M\"{o}bius band $B_2$ whose boundary is disjoint from $\partial B_1$.
Let $A$ be an annulus between $\partial B_1$ and $\partial B_2$ on $\widehat{S}$.
Then $B_1\cup A\cup B_2$ is a Klein bottle $\widehat{F}$ in $K(\alpha)$, which meets the core of $V_\alpha$ in two points
(after a perturbation).
Then $F=\widehat{F}\cap E(K)$ gives a twice-punctured Klein bottle in $E(K)$.
By attaching a suitable annulus on $\partial E(K)$ to $F$ along their boundaries, we obtain a closed non-orientable surface in $S^3$, a contradiction. 

(2) By \cite[Lemma 1.4]{W}, such family contains at most $s/2+2$ edges.
Furthermore, if it contains $s/2+2$ edges, then we can assume that
there are two $S$-cycles with disjoint label pairs $\{1,2\}$ and $\{s/2+1,s/2+2\}$ in the family (see \cite[Figure 1]{W}).
By the same construction as in (1), we have two M\"{o}bius bands $B_1$ and $B_2$ and an annulus $A$ on $\widehat{S}$.
In $\widehat{S}$, two vertices $u_{i}$ and $u_{s-i+3}$ are connected by an edge in the family for $i=3,4,\dots,s/2$.
Hence $\mathrm{Int}\,A$ contains an even number of vertices.
Then $B_1\cup A\cup B_2$ is a Klein bottle which meets $V_\alpha$ in an even number of meridian disks.
This leads to a contradiction as before.

If the family contains $s/2+1$ edges, then it contains an $S$-cycle at one end \cite[Lemma 1.4]{W}.

(3) If $s\ge 4$, then this is \cite[Lemma 2.3(1)]{GW}.
If $s=2$, then any negative edge in $G_T$ has the same label at its endpoints.
Thus if a family of parallel negative edges contains $3$ edges, then there would be two edges which are parallel in both graphs.
This contradicts Lemma \ref{lem:common}.

(4) is \cite[Lemma 1.7]{W}.
\end{proof}

Suppose $q_i=s$.
Let $e_1,e_2,\dots,e_s$ be the edges of $Q_i$, numbered successively.
We may assume that $e_j$ has label $j$ at $v_1$ for $1\le j\le s$.
Then these edges define a permutation $\sigma$ of the set $\{1,2,\dots,s\}$ such that
$e_j$ has label $\sigma(j)$ at $v_2$.
In fact, $\sigma(j)\equiv j+k \pmod{s}$ for some even $k$.
This $\sigma$ is called the \textit{associated permutation\/} to $Q_i$.
By the parity rule, $\sigma$ has at least two orbits, and all orbits have the same length.
According to the orbits of $\sigma$, the edges of $Q_i$ form disjoint cycles on $\widehat{S}$.

\begin{lemma}\label{lem:orbit}
Each of these cycles is essential on $\widehat{S}$.
\end{lemma}

\begin{proof}
This is \cite[Lemma 2.3]{G}.
\end{proof}

%%%%%%%%%%%%%%%%%%%%%%%%%%%%%%%%%%%%%%%%%%%%%%%%%%%%%%%%%%%%%%%%%%%%%%%%%%%%%%%%%%%%%%%%%%%%%%%%%%%%%%%%%%%%%%%%%%%%%%%%%%%%
\section{The case that $s\ge 4$}\label{sec:s4}

In this section, we will show that the case where $s\ge 4$ is impossible.
Recall that $q_1\le s/2+1$ and $q_i\le s$ for $i\ge 2$ by Lemma \ref{lem:GT}.

\begin{lemma}
$q_1=s/2$ or $s/2+1$.
\end{lemma}

\begin{proof}
Since $2q_1+q_2+q_3+q_4+q_5=5s$, $q_1\ge s/2$.
The result follows from this.
\end{proof}

We distinguish two cases.

%%%
\subsection{$q_1=s/2$}

\begin{figure}[tb]
%\blankbox{1.0\columnwidth}{2in}
\includegraphics*[scale=0.4]{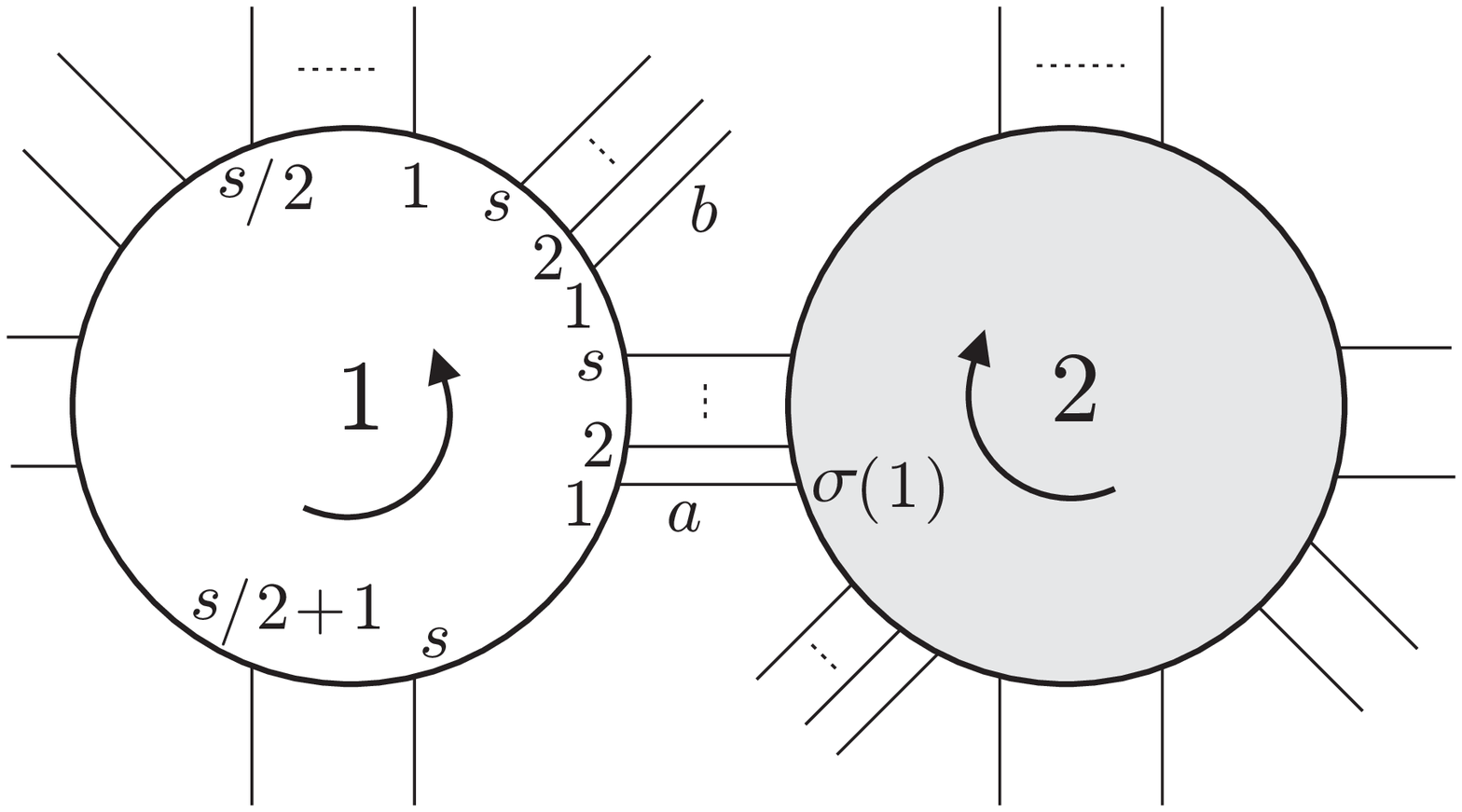}
\caption{}\label{fig:s4case1}
\end{figure}

Then $G_T\cong G(s/2,s,s,s,s)$.
We can assume that $G_T$ has labels as in Figure \ref{fig:s4case1}.

Let $\sigma$ be the associated permutation to $Q_2$ such that an edge in $Q_2$ has label $i$ at $v_1$ and label $\sigma(i)$ at $v_2$.
Notice that the associated permutations to $Q_3,Q_4$ and $Q_5$ are equal to $\sigma$.

\begin{lemma}\label{lem:case1id}
$\sigma$ is not the identity.
\end{lemma}

\begin{figure}[tb]
%\blankbox{1.0\columnwidth}{2in}
\includegraphics*[scale=0.5]{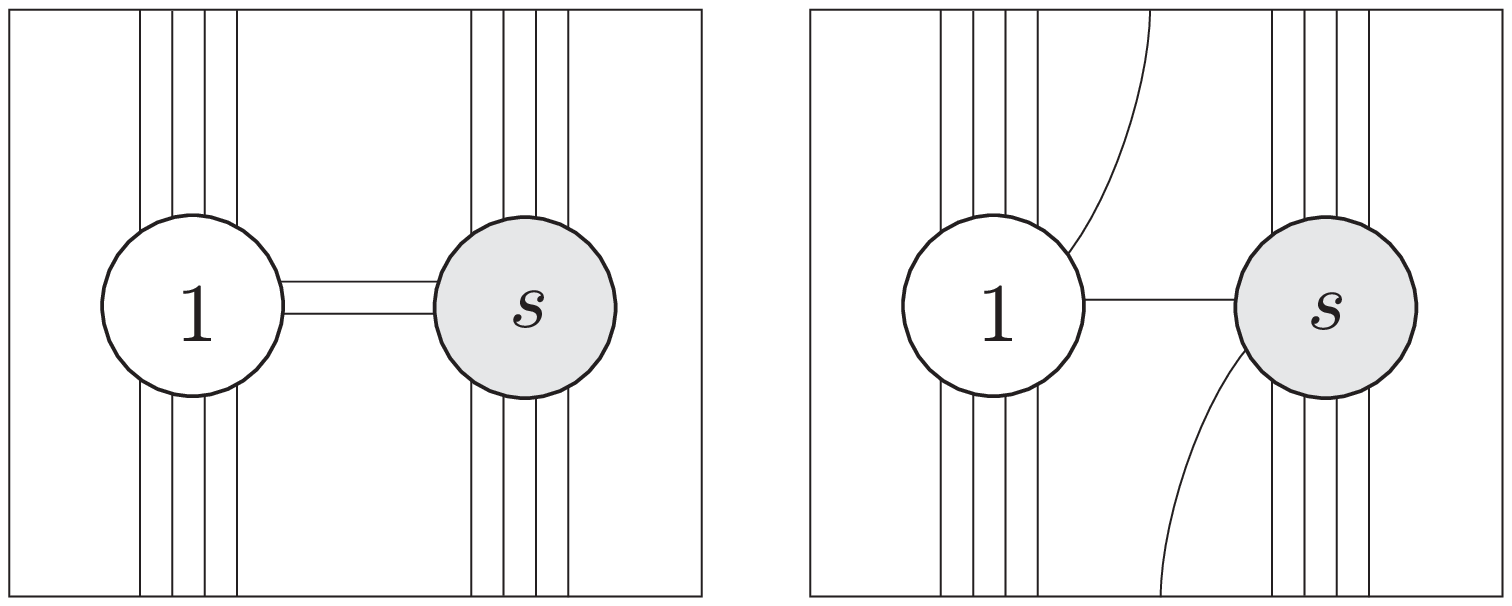}
\caption{}\label{fig:case1id}
\end{figure}

\begin{proof}
Suppose that $\sigma$ is the identity.
Then each vertex of $G_S$ is incident to $4$ loops.
Let $G(1,s)$ be the subgraph of $G_S$ spanned by $u_1$ and $u_s$.
Since $s\ge 4$, $G(1,s)$ has an annulus support on $\widehat{S}$, and
there are only two possibilities for it as in Figure \ref{fig:case1id}, where
the top and bottom edges are identified to form an annuls.

Let $a$ be the $\{1,1\}$-edge in $Q_2$, and let $a_i$ be its endpoint at $v_i$ for $i=1,2$.
Let $e$ (resp.~$f$) be the $\{1,s\}$-loop at $v_1$ ($v_2$), and let $e_1$ and $f_1$ be their endpoints with label $1$. 
Around $v_1$, $a_1$ and $e_1$ are not successive among $5$ occurrences of label $1$, but $a_2$ and $f_1$ are successive among $5$ occurrences of
label $1$ around $v_2$.
By Lemma \ref{lem:jumping}, $a_1$ and $e_1$ are successive among $5$ occurrences of label $1$ around $u_1$, but $a_2$ and $f_1$ are not
successive among $5$ occurrences of label $2$ around $u_1$.
But this is not satisfied in both configurations of Figure \ref{fig:case1id}.
\end{proof}

\begin{lemma}\label{lem:involution}
$\sigma^2$ is the identity.  In particular, each orbit of $\sigma$ has length two, and $\sigma(i)=i+s/2$.
\end{lemma}

\begin{figure}[tb]
%\blankbox{1.0\columnwidth}{2in}
\includegraphics*[scale=0.5]{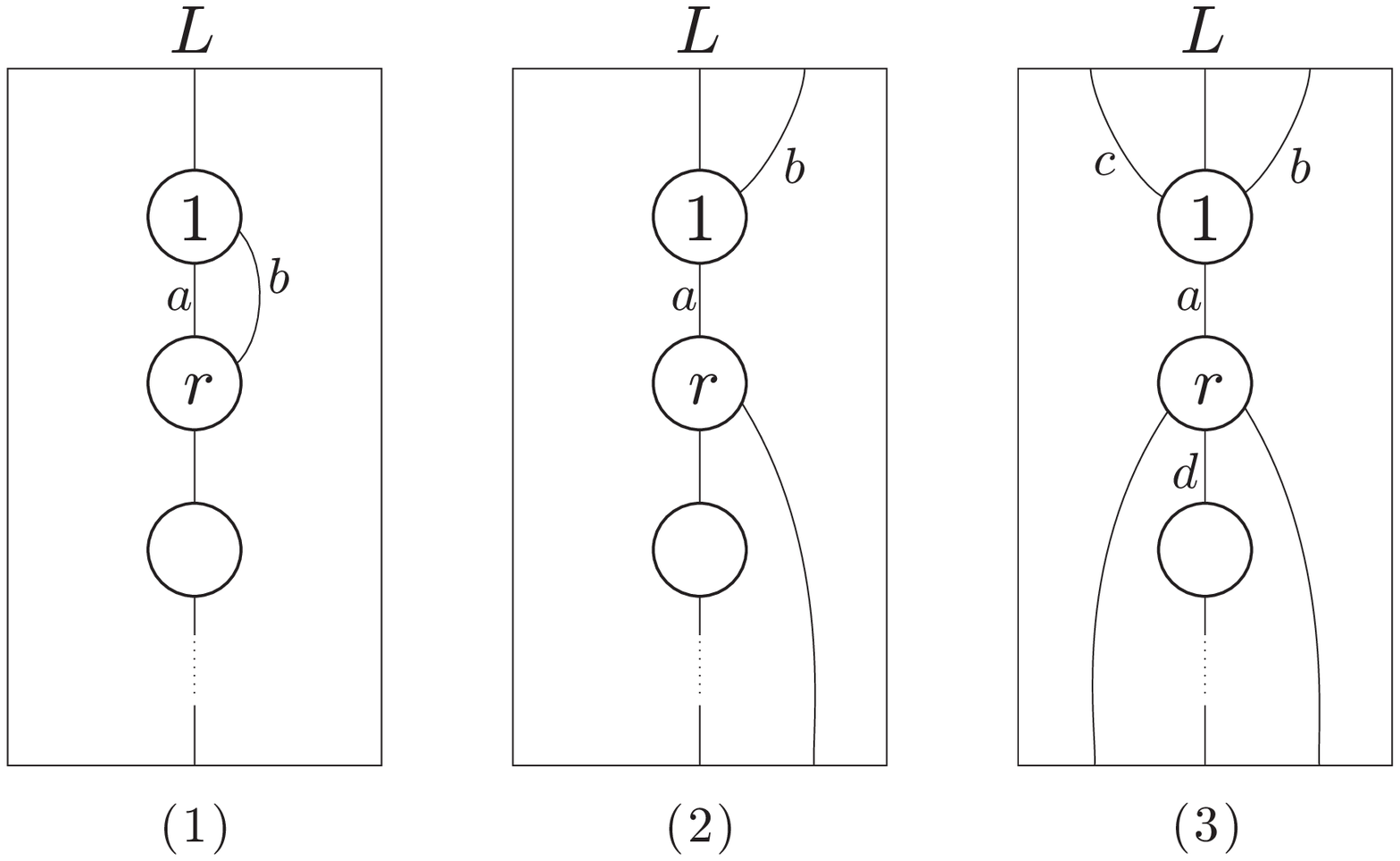}
\caption{}\label{fig:case1inv}
\end{figure}

\begin{proof}
The edges of $Q_2$ form disjoint cycles in $G_S$ according to the orbits of $\sigma$, and such cycle is essential on $\widehat{S}$ by Lemma \ref{lem:orbit}.
Recall that there are at least two such cycles.
Let $L$ be the cycle corresponding to the orbit of $\sigma$ containing $1$.
Let $a$ and $b$ be the edges of $Q_2$ and $Q_3$ with label $1$ at $v_1$, respectively
(see Figure \ref{fig:s4case1}).
Then $L\cup b$ has an annulus support on $\widehat{S}$.
Thus there are two possibilities for $L\cup b$ as in Figure \ref{fig:case1inv}(1) and (2), where $r=\sigma(1)$.
Note that $a$ and $b$ have label $1$ at $u_1$.

For Figure \ref{fig:case1inv}(1), there is another edge $e$ between $a$ and $b$.
Then $e$ is a negative $\{1,r\}$-edge in $G_T$ with label $r$ at $v_1$ and label $1$ at $v_2$.
Although $e$ need not be in $Q_2$, this implies $\sigma(r)=1$, because any family of negative parallel edges
has the same permutation $\sigma$.
Hence $\sigma^2$ is the identity.

For Figure \ref{fig:case1inv}(2), let $c$ and $d$ be the edges of $Q_4$ and $Q_5$ with label $1$ at $v_1$, respectively.
Then there would be a pair of parallel edges among $a$, $b$, $c$, $d$ on $\widehat{S}$.
Thus the above argument gives the conclusion.
\end{proof}

\begin{lemma}\label{lem:half}
The case $q_1=s/2$ is impossible.
\end{lemma}

\begin{proof}
By Lemma \ref{lem:involution}, $G_S^+$ consists of $s/2$ components, each of which has an annulus support.
In fact, any component consists of two vertices and $8$ edges.
Thus there are at least $4$ mutually parallel positive edges.
But this is impossible by Lemma \ref{lem:GS}(1).
\end{proof}

%%%
\subsection{$q_1=s/2+1$}

Since $q_2+q_3+q_4+q_5=4s-2$, at least two of $q_i$ are $s$.
By the parity rule, $q_2+q_3$ and $q_4+q_5$ are even.
Thus we may assume that $G_T\cong G(s/2+1,s,s,s,s-2)$ or $G(s/2+1,s,s,s-1,s-1)$ without loss of generality.

\begin{lemma}
$G(s/2+1,s,s,s-1,s-1)$ is impossible.
\end{lemma}

\begin{proof}
Assume $G_T\cong G(s/2+1,s,s,s-1,s-1)$.
By the parity rule, each edge of $Q_2$ has labels of the same parity.
But then each edge of $Q_4$ has labels of opposite parities at its ends.
This contradicts the parity rule.
\end{proof}

Thus $G_T\cong G(s/2+1,s,s,s,s-2)$.
We can assume that the edges of $Q_2$ have labels $1,2,\dots,s$ at $v_1$ as in Figure \ref{fig:s4case1}.
Then there is an $S$-cycle with label pair $\{s/2,s/2+1\}$ among loops at $v_1$.
Let $\sigma$ be the associated permutation to $Q_2$ as before.
Notice that $Q_3$ and $Q_5$ also associate to $\sigma$.
(Although $q_5=s-2$, the associated permutation is obviously defined.)

\begin{lemma}\label{lem:id}
$\sigma$ is the identity.
\end{lemma}

\begin{proof}
Among loops at $v_2$, there is an $S$-cycle with label pair $\{r+s/2-1,r+s/2\}$, where $r=\sigma(1)$.
Assume that $\sigma$ is not the identity.  Then $r$ is odd by the parity rule, and $r\ge 3$.

If $s=4$, then $r=3$.  Hence $G_T$ contains two $S$-cycles with label pairs $\{2,3\}$ and $\{4,1\}$.
But this contradicts Lemma \ref{lem:GT}(1).
Hence $s>4$.

We claim that $\sigma^2$ is the identity.
To see this, assume not.
As in the proof of Lemma \ref{lem:involution}, the edges of $Q_2$ form disjoint cycles according to the orbits of $\sigma$.
Let $L$ be the cycle through $u_1$.
Notice that $L$ contains at least three vertices.
Let $a$ and $b$ be the edges of $Q_2$ and $Q_3$ with label $1$ at $v_1$.
Then $L\cup b$ has an annulus support.
If $b$ is parallel to $a$, then $\sigma^2$ would be the identity.
Hence $L\cup b$ is as shown in Figure \ref{fig:case1inv}(2).
Let $c$ be the edge of $Q_5$ with label $1$ at $v_1$.
If $c$ is parallel to $a$ or $b$ on $\widehat{S}$, then $\sigma^2$ would be the identity again.
Hence $c$ is located as in Figure \ref{fig:case1inv}(3).
Let $d$ and $e$ be the edges of $Q_2$ and $Q_3$ with label $r$ at $v_1$.
Of course, $d$ is in $L$.
But $e$ must be parallel to $d$.
Hence $\sigma^2$ would be the identity, a contradiction.

Therefore $r=s/2+1$.
Since $r$ is odd, $s\equiv 0\pmod{4}$.
The edges of $Q_2$ form cycles of length two on $\widehat{S}$, and
there are at least $4$ such cycles.
In particular, $u_1$ and $u_{s/2+1}$ lie on the same cycle, and so do $u_{s/2}$ and $u_{s}$.

Also, $G_T$ contain two $S$-cycles with label pairs $\{s/2,s/2+1\}$ and $\{s,1\}$.
But we cannot locate the edges of these $S$-cycles to form essential cycles on $\widehat{S}$ simultaneously.
This contradicts Lemma \ref{lem:common}(3).
\end{proof}

\begin{lemma}
$s=4$.
\end{lemma}

\begin{proof}
Assume $s>4$.
Notice that the associated permutation $\tau$ to $Q_4$ is defined as $\tau(i)=i-2$.
Hence the edges of $Q_4$ form two cycles on $\widehat{S}$, each of which contains at least three vertices.
By Lemma \ref{lem:id}, each vertex of $G_S$ is incident to a loop.
Hence there would be a trivial loop.
\end{proof}

\begin{lemma}
The case $q_1=s/2+1$ is impossible.
\end{lemma}

\begin{proof}
In $G_S$, $u_1$ and $u_4$ are incident to $3$ loops respectively, and $u_2$ and $u_3$ are incident to two loops.
Also, $G_T$ contains two $S$-cycles with label pair $\{2,3\}$.
The edges of them give $4$ edges between $u_2$ and $u_3$ in $G_S$.
There are two edges between $u_1$ and $u_3$, and $u_2$ and $u_4$, which belong to $Q_4$ in $G_T$.
Hence they are not parallel in $G_S$ by Lemma \ref{lem:common}(1).
Thus $G_S$ is as shown in Figure \ref{fig:s4}.

\begin{figure}[tb]
%\blankbox{1.0\columnwidth}{2in}
\includegraphics*[scale=0.5]{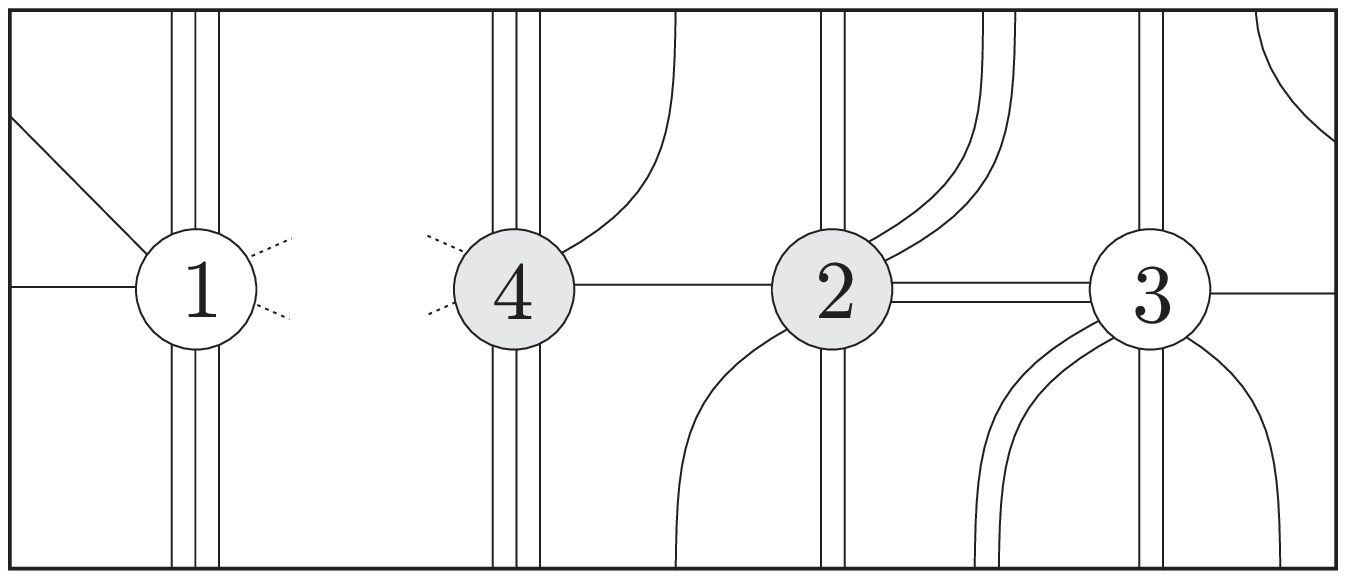}
\caption{}\label{fig:s4}
\end{figure}

Let $a,b,d,e$ be the edges of $Q_2,Q_3,Q_4,Q_5$ with label $1$ at $v_1$.
Also, let $c$ be the loop at $v_1$ with label $1$.
By Lemma \ref{lem:jumping}, these edges appear in the order $a,c,e,b,d$ around $u_1$ in some direction.
Let $f$ be the edge of $Q_4$ with label $3$ at $v_1$.
Then the endpoints of $b$ and $f$ at $v_2$ are consecutive among $5$ occurrences of label $1$.
Also, their endpoints at $u_1$ are consecutive among $5$ occurrences of label $2$.
This contradicts Lemma \ref{lem:jumping}.
\end{proof}

%%%%%%%%%%%%%%%%%%%%%%%%%%%%%%%%%%%%%%%%%%%%%%%%%%
\section{The case that $s=2$}\label{sec:s2}

In this section, we show that there is only one possible pair of the graphs, which will lead to the determination of knots in the next section.

The reduced graphs $\overline{G}_S$ and $\overline{G}_T$ are subgraphs of the graph as shown in Figure \ref{fig:t2}.
We use $p_i$ (resp.~$q_i$) to denote the weight of edge in $\overline{G}_S$ (resp.~$\overline{G}_T$).
In Figure \ref{fig:t2}, $q_i$ is indicated, but $p_i$ is assigned at the same place as $q_i$.
Also, we say $G_S\cong G(p_1,p_2,p_3,p_4,p_5)$, similarly to $G_T$.
Recall that $p_1\le 3$ and $p_i\le 2$ for $i=2,3,4,5$ by Lemma \ref{lem:GS}.
Since $2p_1+p_2+p_3+p_4+p_5=10$, $p_1\ge 1$.

\begin{lemma}
$p_1=1$ is impossible.
\end{lemma}

\begin{proof}
If $p_1=1$, then $p_i=2$ for any $i\ne 1$.
Then $G_T\cong G(4,2,0,0,0)$ or $G(4,1,1,0,0)$.
But these are eliminated as in the proof of Lemma \ref{lem:case1id}.
\end{proof}

Thus $p_1=2$ or $3$.

\begin{figure}[tb]
%\blankbox{1.0\columnwidth}{2in}
\includegraphics*[scale=0.5]{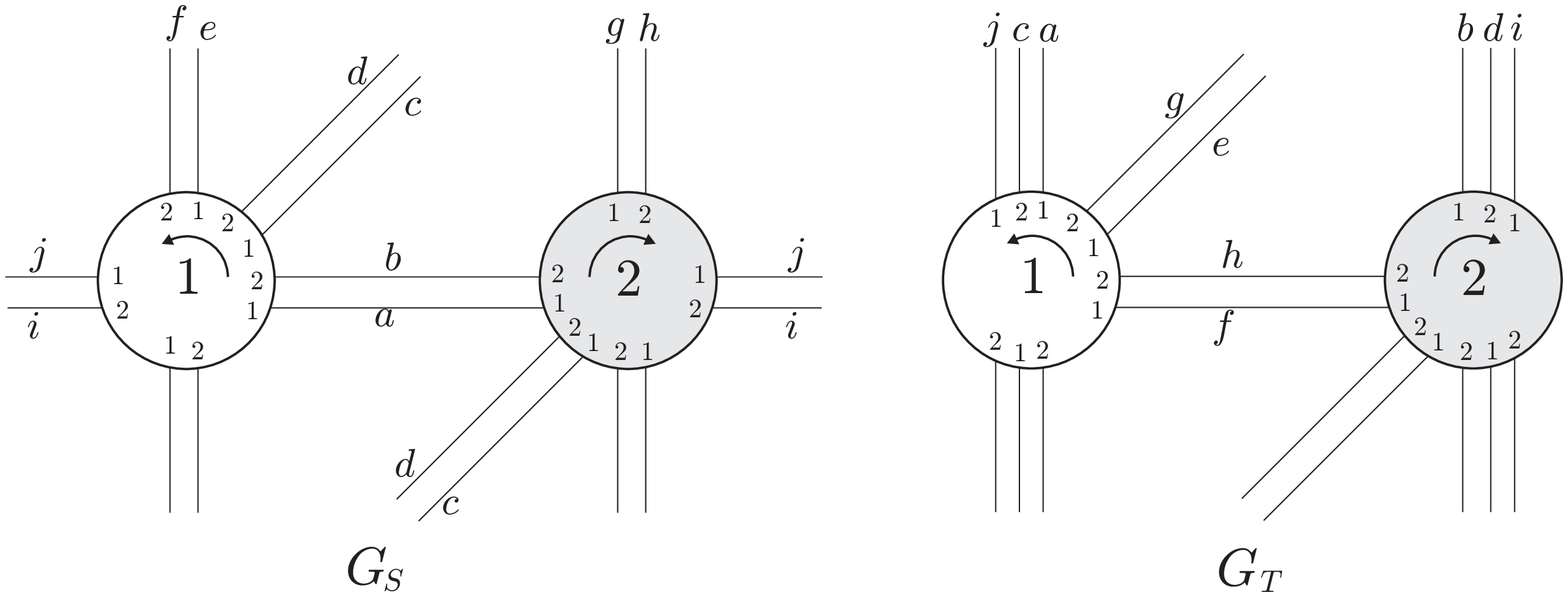}
\caption{}\label{fig:s2t2final}
\end{figure}

\begin{lemma}\label{lem:p2}
If $p_1=2$, then the graphs are as shown in Figure \ref{fig:s2t2final}.
The correspondence between the edges in $G_S$ and $G_T$ is uniquely determined up to symmetry.
\end{lemma}

\begin{figure}[tb]
%\blankbox{1.0\columnwidth}{2in}
\includegraphics*[scale=0.6]{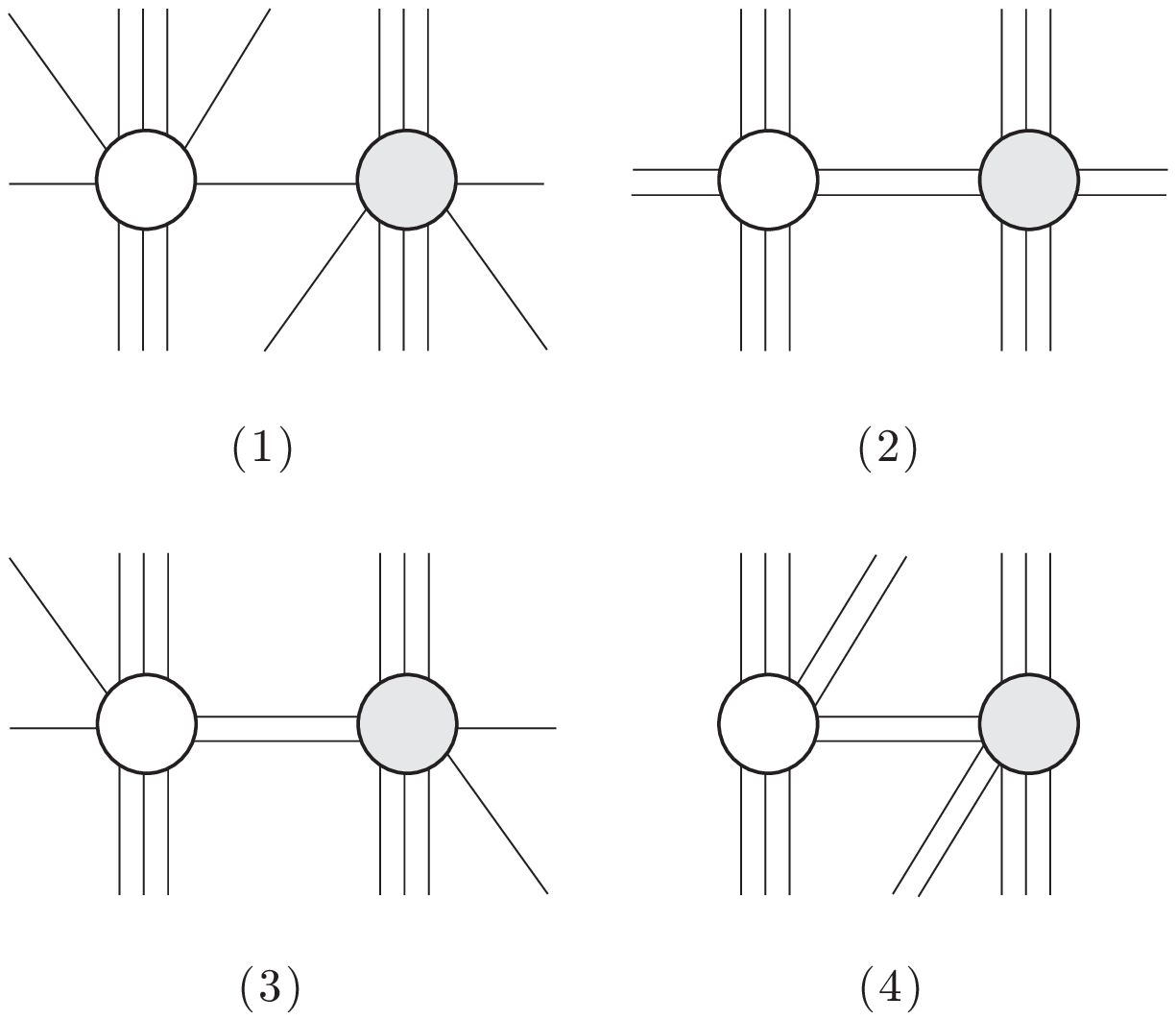}
\caption{}\label{fig:s2t2case2}
\end{figure}

\begin{figure}[tb]
%\blankbox{1.0\columnwidth}{2in}
\includegraphics*[scale=0.7]{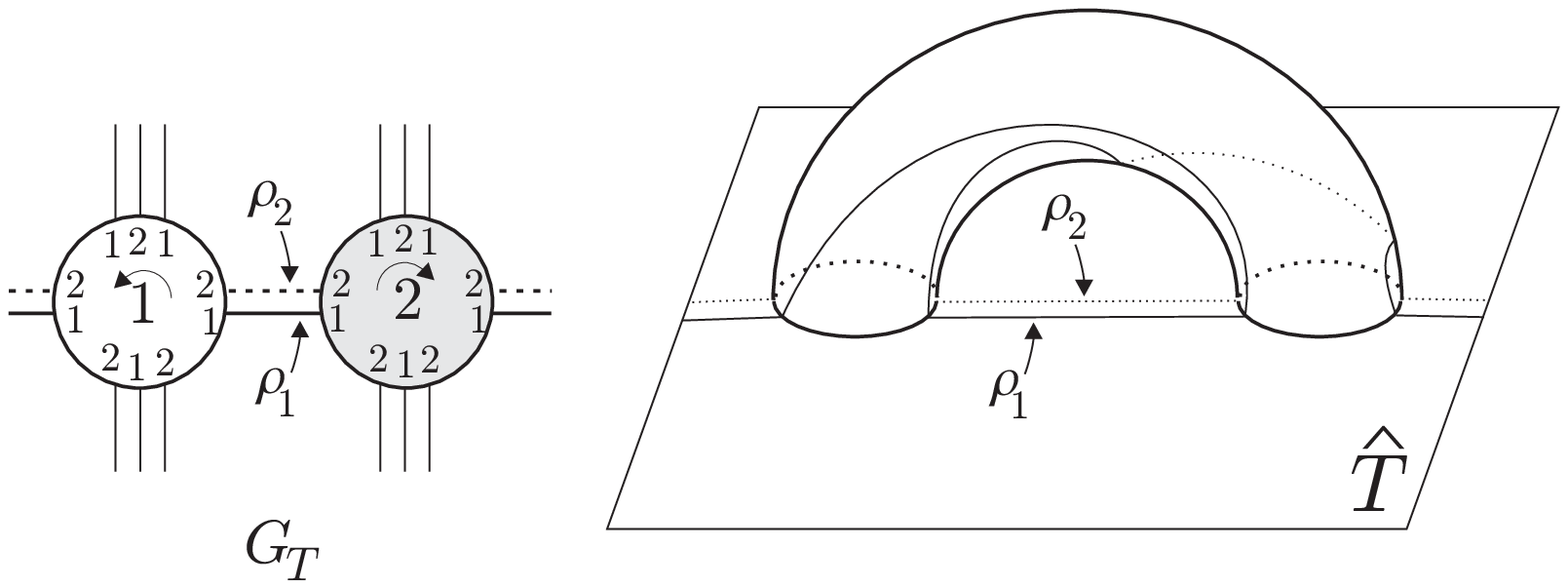}
\caption{}\label{fig:s2t2case2irr}
\end{figure}

\begin{proof}
Suppose $p_1=2$.
By the parity rule, $p_2+p_3$ and $p_4+p_5$ is even.
Hence we can assume that $p_2+p_3=4$ and $p_4+p_5=2$.
Then $p_2=p_3=2$ by Lemma \ref{lem:GS}.
If $p_4=p_5=1$, then the labels in $G_S$ contradicts the parity rule.
Thus we can assume that $G_S\cong G(2,2,2,2,0)$, and then 
there are $4$ possibilities for $G_T$ up to symmetry as in Figure \ref{fig:s2t2case2}.

If $G_T$ is (1), then the endpoints of two negative $\{1,1\}$-edges at $v_1$ are not successive among
$5$ occurrences of label $1$.
But these points are also not successive at $u_1$.
This contradicts Lemma \ref{lem:jumping}.

To eliminate (2), notice that $G_S$ contains two $S$-cycles $\rho_1$ and $\rho_2$ whose faces lie on the same side of $\widehat{T}$.
From the labeling of $G_S$, we can determine the edges of $\rho_i$ in $G_T$ as in Figure \ref{fig:s2t2case2irr}.
Then it is impossible to connect these edges of $\rho_1$ and $\rho_2$ on $\partial V_\beta$ simultaneously.

Clearly, (3) contradicts the parity rule.

For (4), the edge correspondence is determined in the same way as \cite[Lemma 5.3]{GW} by using the fact that the jumping number is two.
Notice that this pair coincides with that of \cite[Figures 5.2(e) and 5.3(b)]{GW}, when $G_T$ is regarded as a graph in an annulus.
\end{proof}

\begin{lemma}\label{lem:p3}
If $p_1=3$, then the graphs are the same as in Figure \ref{fig:s2t2final} with exchanging $G_S$ and $G_T$.
\end{lemma}

\begin{proof}
We may assume that $(p_2+p_3,p_4+p_5)=(2,2)$ or $(4,0)$ by symmetry.
In the former case, there are three possibilities for $G_S$, up to equivalence, as shown in Figure \ref{fig:s2t2case2}(1), (2) and (3).

(3) is impossible by the parity rule.
If $G_S$ is (2), then the labeling of $G_S$ implies that $G_T$ contains an $S$-cycle at each vertex.
It is easy to see that their faces lie on the same side of $\widehat{S}$.
In fact, $G_T\cong G(2,2,2,2,0)$.
Hence the argument in the proof of Lemma \ref{lem:p2} works again (with an exchange of roles between $G_S$ and $G_T$).
If $G_S$ is (1), then the endpoints of two $\{1,1\}$ negative edges at $u_1$ are not successive among $5$ occurrences of label $1$.
Hence Lemma \ref{lem:jumping} implies that $(q_2+q_3,q_4+q_5)=(6,0)$, up to symmetry. 
But this is impossible since $q_i\le 2$ for $i\ne 1$.

Thus $(p_2+p_3,p_4+p_5)=(4,0)$, giving $p_2=p_3=2$.
Hence $q_1=2$, and so we can assume that $(q_2+q_3,q_4+q_5)=(6,0)$ or $(4,2)$.
But $(6,0)$ is impossible again.
By using the parity rule, it is easy to see that $G_T$ is as in Figure \ref{fig:s2t2final} (with exchanging $G_S$ and $G_T$).
\end{proof}

We have thus shown that there is only one possibility for the pair $\{G_S,G_T\}$ realizing the distance $5$, as
shown in Figure \ref{fig:s2t2final} (or with an exchange of $G_S$ and $G_T$).

%%%%%%%%%%%%%%%%%%%%%%%%%%%%%%%%%%%%%%%%%%%%%%%%%%%%%%%%%%%%%%%%%%%%%%%%%%%%%%%%%%%%%%%%%%%%%%%%%%%%%%%%%%%%%%%%%%%%%%%%%%
\section{Determining the knot}\label{sec:knot}

In this section, we will show that the only possibility for $\{G_S,G_T\}$ corresponds to the Eudave-Mu\~{n}oz knot $k(2,-1,n,0)$.
We follow the construction in \cite[Section 5]{GL3}, and use the same notation as far as possible.
We assume that the pair is as in Figure \ref{fig:s2t2final} with exchanging $G_S$ and $G_T$ to fit the notation to \cite{GL3}.
As shown in Lemmas \ref{lem:p2} and \ref{lem:p3}, $G_S$ and $G_T$ can be exchanged.
Hence we cannot say which of $\alpha$ and $\beta$ is non-integral with respect to
the original framing of $K$ in the following argument.

\begin{figure}[tb]
%\blankbox{1.0\columnwidth}{2in}
\includegraphics*[scale=0.5]{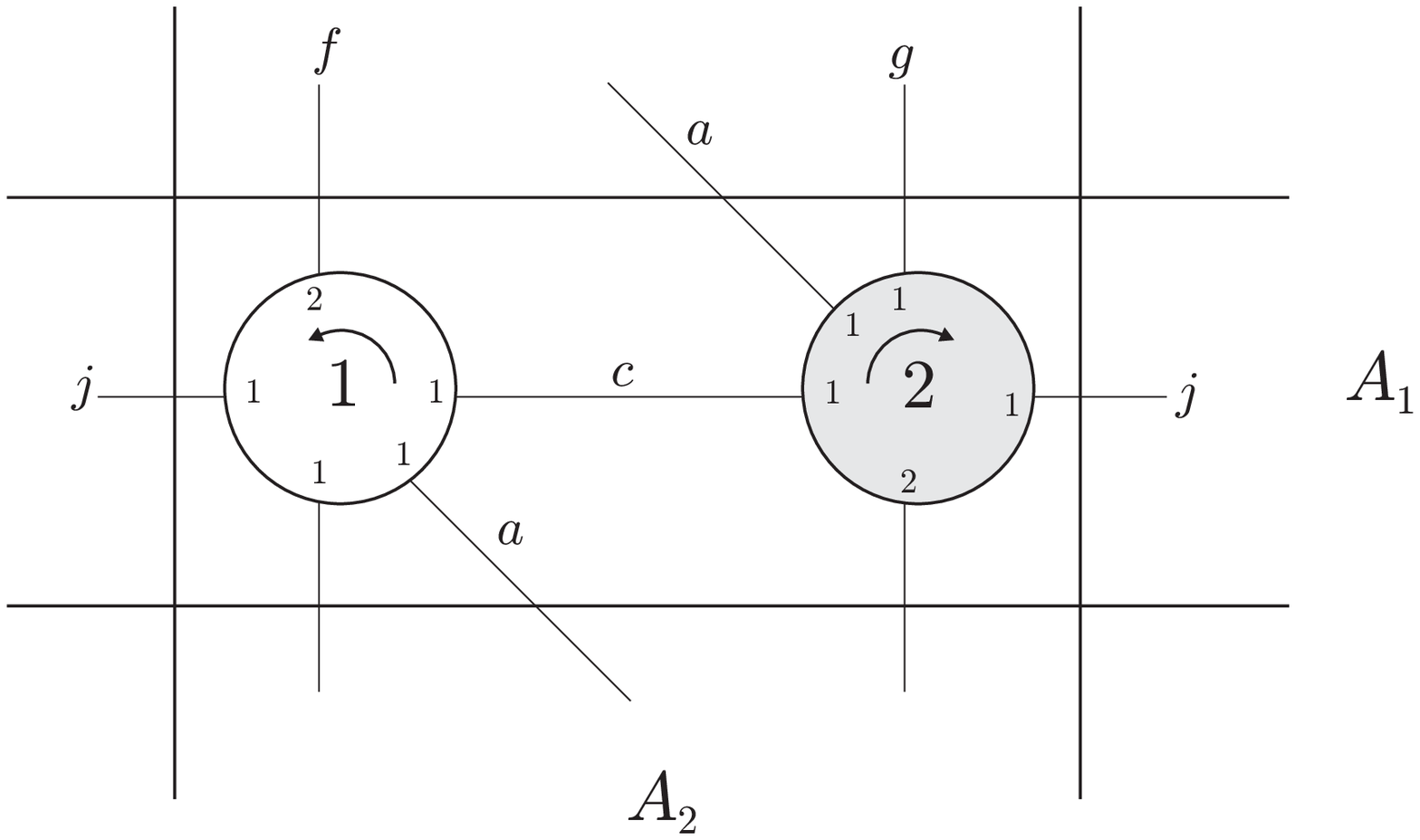}
\caption{}\label{fig:asupp}
\end{figure}

Let $f_1$ and $f_2$ be the bigons in $G_S$ bounded by the edges $c$ and $j$, $a$ and $c$, respectively.
Also, let $f_3$ be the $3$-gon bounded by $a,g$ and $f$.
Up to homeomorphism, we can assume that these edges are as in Figure \ref{fig:asupp} on $\widehat{T}$.
For $i=1,2$,
let $A_i$ be an annulus support of the edges of $f_i$, and let $A_i'=\mathrm{cl}(\widehat{T}-A_i)$.

\begin{figure}[tb]
%\blankbox{1.0\columnwidth}{2in}
\includegraphics*[scale=0.7]{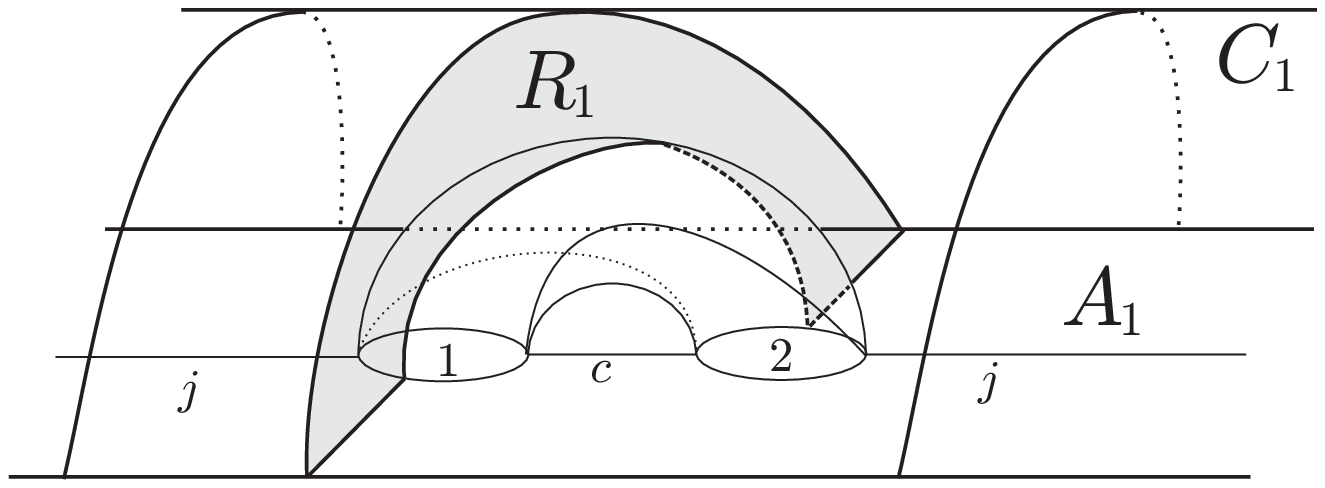}
\caption{}\label{fig:r1}
\end{figure}

\begin{figure}[tb]
%\blankbox{1.0\columnwidth}{2in}
\includegraphics*[scale=0.7]{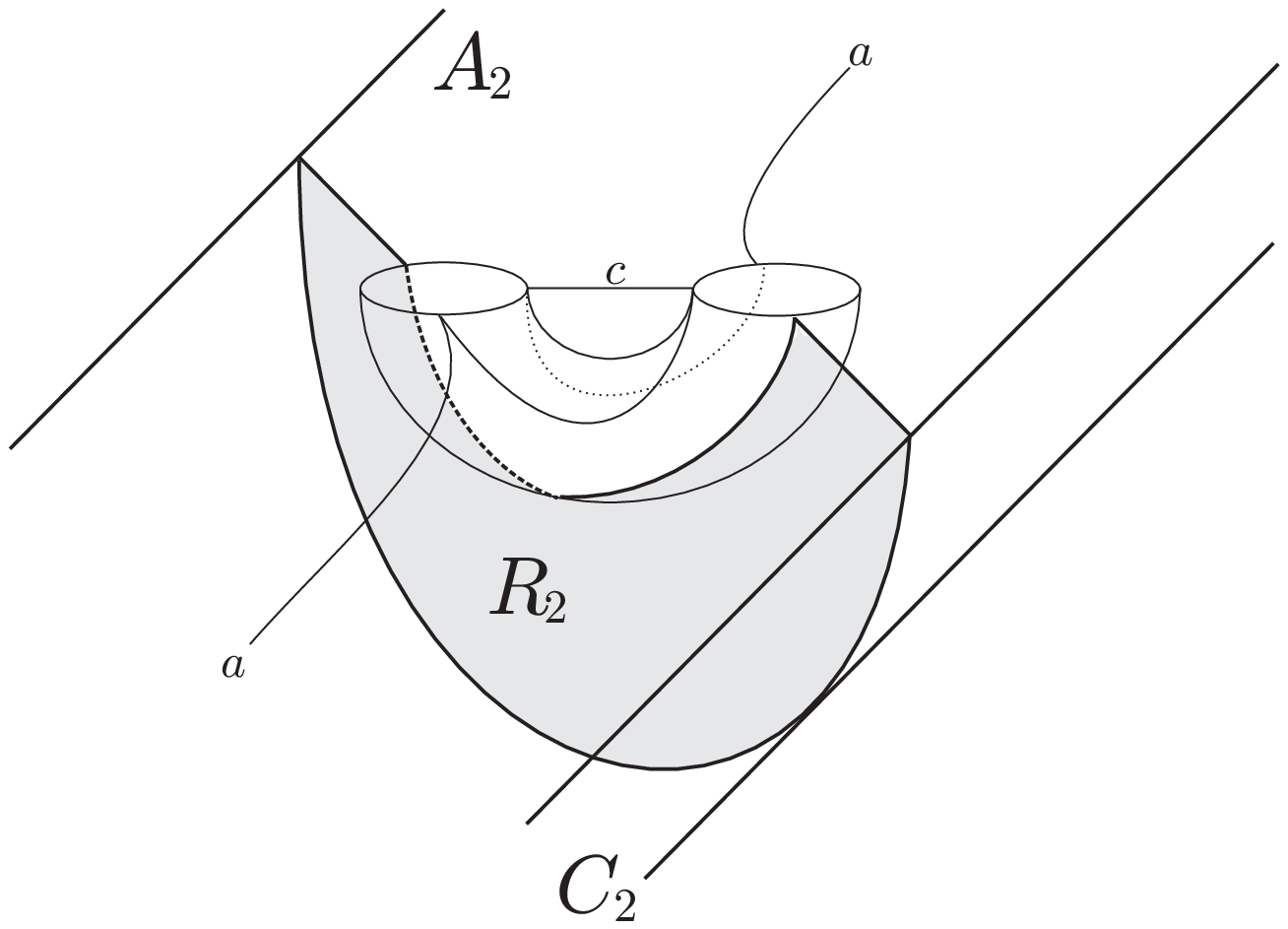}
\caption{}\label{fig:r2}
\end{figure}

Let $K(\beta)=M_1\cup_{\widehat{T}}M_2$, and let $H_i=V_\beta\cap M_i$.
We may assume that $f_i$ lies in $M_i$.
Then $W_i=N(A_i\cup H_i\cup f_i) \subset M_i$ is a solid torus \cite[Lemma 3.7]{GL}, and moreover $M_i=W_i\cup W_i'$, where
$W_i'$ is a solid torus, is a Seifert fibered manifold over the disk with two exceptional fibers, one of which has index two \cite[Lemma 3.8]{GL}.
Let $\partial W_i=A_i\cup C_i$.
Take a disk (rectangle) $R_i$ in $\mathrm{cl}(W_i-H_i)$ as in \cite{GL3}.
See Figures \ref{fig:r1} and \ref{fig:r2}.

Now, regard $\widehat{T}\cup C_1\cup C_2\cup R_1\cup R_2 \cup V_\beta$ as a subset of $S^3$ as follows.
Embed $\widehat{T}$ as a standard torus in $S^3$, which splits $S^3$ into two solid tori $U_1$ and $U_2$.
We assume that the components of $\partial A_1$ are longitudes of $U_1$, and those of $\partial A_2$ are meridians of $U_1$.
Then $C_i$ is identified with the obvious annulus in $U_i$, separating $U_i$ into two solid tori $V_i$ and $V_i'$, where
$\partial V_i=A_i\cup C_i$ and $\partial V_i'=A_i'\cup C_i$.
Finally, embed $V_\beta$ in the obvious way, and $R_i$ as in Figures \ref{fig:r1} and \ref{fig:r2}.

\begin{figure}[tb]
%\blankbox{1.0\columnwidth}{2in}
\includegraphics*[scale=0.5]{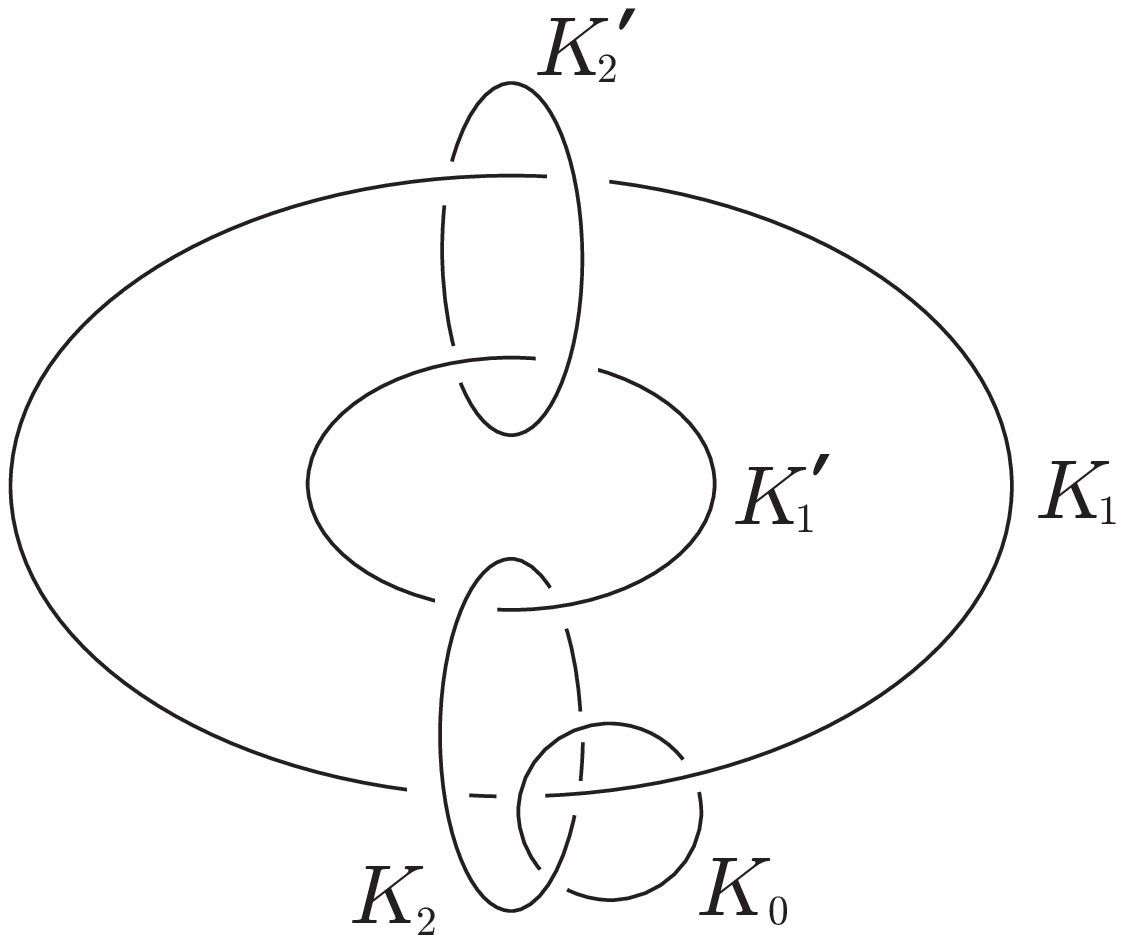}
\caption{}\label{fig:link}
\end{figure}

Let $K_i$ be a core of $V_i$ disjoint from $H_i\cup R_i$, and let $K_i'$ be a core of $V_i'$.
Notice that $N(\widehat{T}\cup C_1\cup C_2)$ is the exterior of the link $L_0=K_1\cup K_1'\cup K_2\cup K_2'$.
Since $W_i'$ is a solid torus, it is obtained from $V_i'$ by some Dehn surgery on $K_i'$.
Similarly, $(W_i; H_i,R_i)$ is obtained from $(V_i;H_i,R_i)$ by some Dehn surgery on $K_i$.
Let $K_0$ be a core of $V_\beta \subset S^3$.
Thus $(K(\beta),K_\beta)$ is obtained from $(S^3,K_0)$ by some Dehn surgery on $L_0$.
In other words, $E(K)$ is obtained from $E(K_0)$ by Dehn surgery on $L_0$.
See Figure \ref{fig:link}.

\begin{lemma}
In the surgery description of $E(K)$ in Figure \ref{fig:link},
$\alpha=-3/5$ and $\beta=1/0$.
\end{lemma}

\begin{figure}[tb]
%\blankbox{1.0\columnwidth}{2in}
\includegraphics*[scale=0.55]{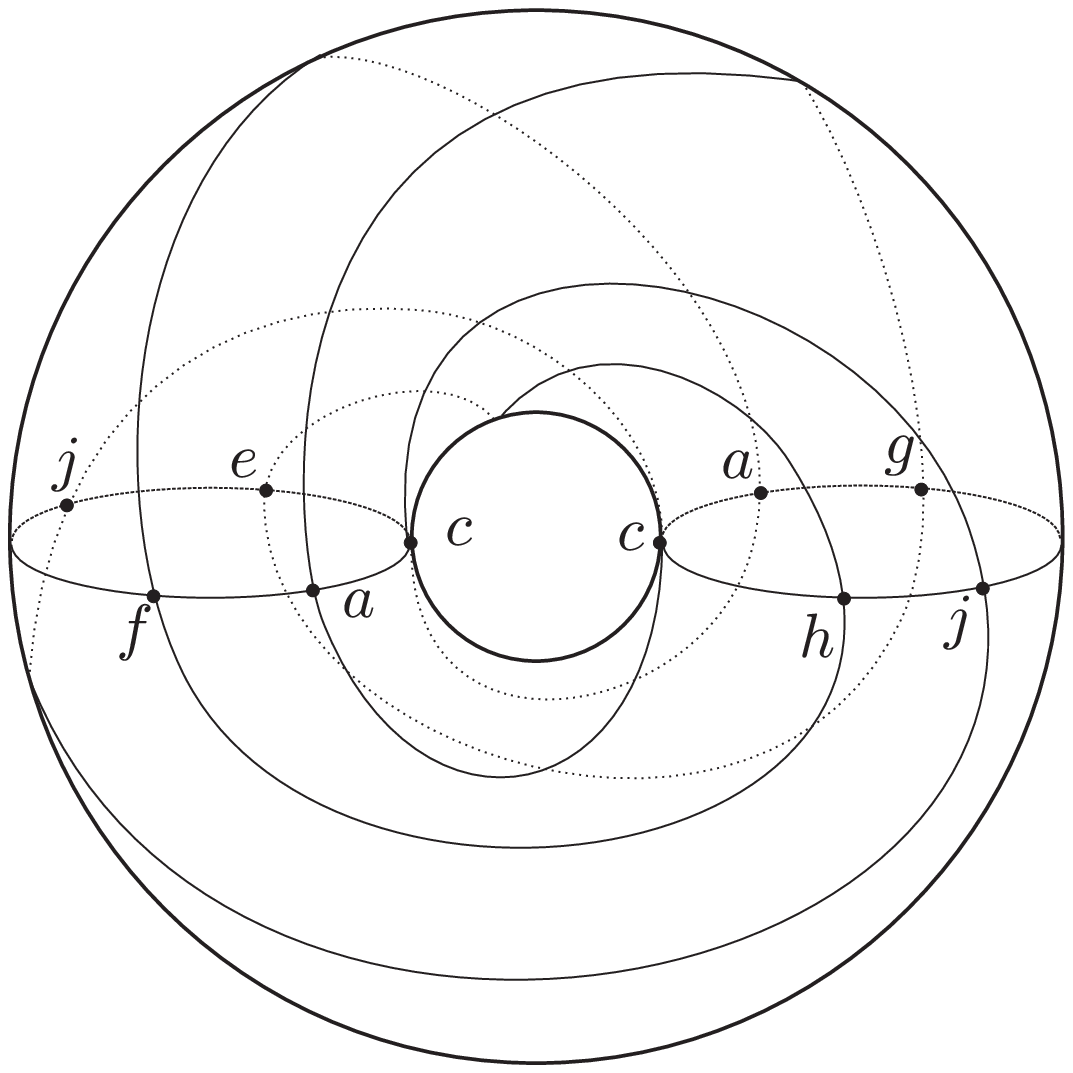}
\caption{}\label{fig:beta}
\end{figure}

\begin{proof}
This is obvious for $\beta$.
For $\alpha$, we draw $\partial u_1$ on $V_\beta$.
See Figure \ref{fig:beta}.
Thus $\alpha=-3/5$.
\end{proof}

\begin{lemma}
The filling slopes for $(K_1,K_1',K_2,K_2')$ are $(-2,3,-2,p/q)$, where $p/q\ne 0/1$.
\end{lemma}

\begin{figure}[tb]
%\blankbox{1.0\columnwidth}{2in}
\includegraphics*[scale=0.9]{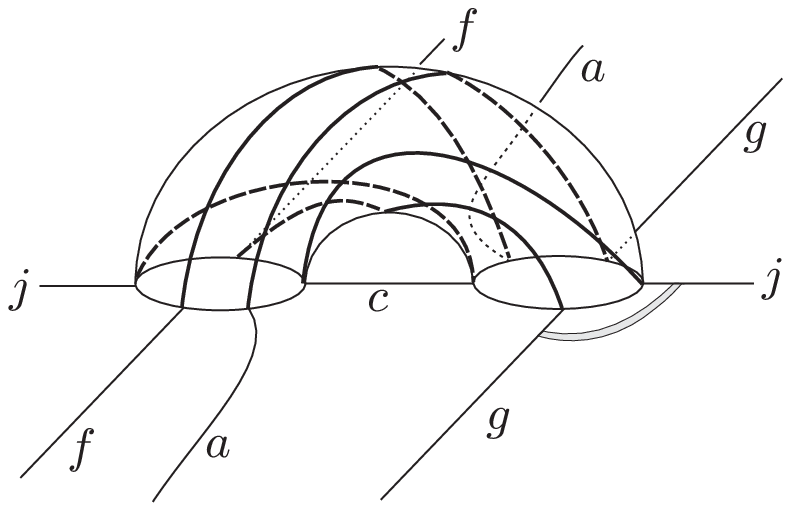}
\caption{}\label{fig:f3}
\end{figure}

\begin{proof}
In Figure \ref{fig:r1}, it is easy to see that there is an annulus between $\partial f_1$ and the curve of slope $-2/1$ on $\partial N(K_1)$.
This is similar for $K_2$.
Although we cannot determine the filling slope on $K_2'$, it cannot be $0/1$, since a regular fiber of the Seifert fibration in $M_2$
is isotopic to a core of $A_2$ (see \cite[Lemma 3.8]{GL}).

To determine the slope on $K_1'$, we need the $3$-gon $f_3$.
Notice that $\partial f_3$ is non-separating on $\partial W_1'$.
Hence $f_3$ gives a meridian disk of $W_1'$.
Let $f$ be the disk obtained by band summing $f_3$ and $f_1$ along the band on $\widehat{T}$ as shown in Figure \ref{fig:f3}.
Then it is easy to see that there is an annulus between $\partial f$ and the curve of slope $3/1$ on $\partial N(K_1')$.
\end{proof}

\begin{figure}[tb]
%\blankbox{1.0\columnwidth}{2in}
\includegraphics*[scale=0.5]{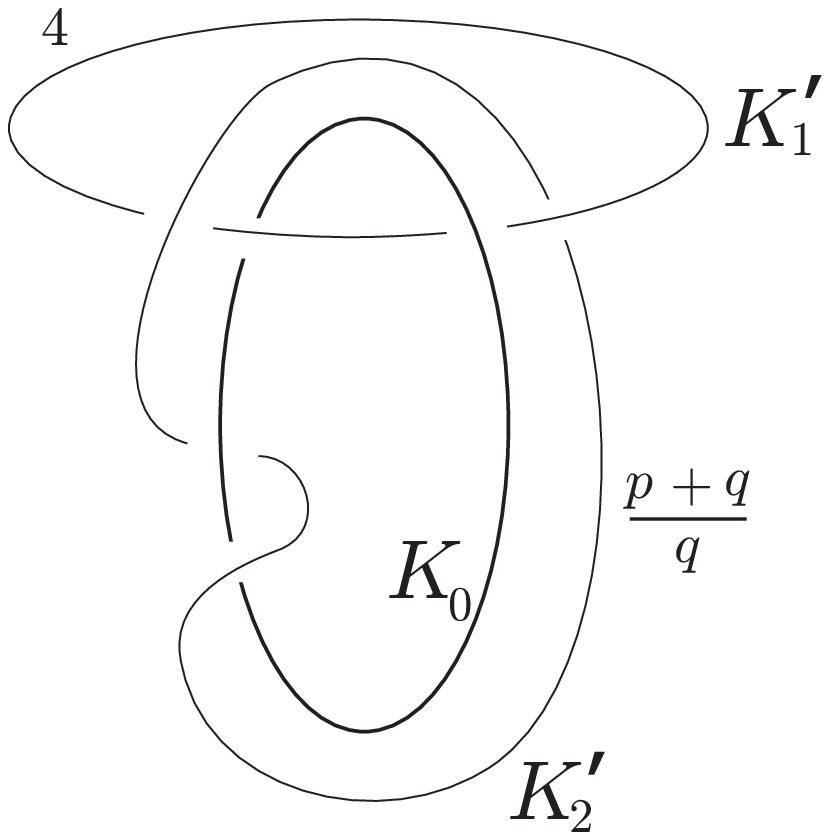}
\caption{}\label{fig:surgery}
\end{figure}

\begin{lemma}\label{lem:surgery}
$E(K)$ has the surgery description as in Figure \ref{fig:surgery}, where
the filling slopes for $(K_1',K_2')$ are $(4,(p+q)/q)$.
Also, $\alpha=1/2$ and $\beta=3$ on $K_0$.
\end{lemma}

\begin{proof}
In the surgery description in Figure \ref{fig:link},
perform a Kirby-Rolfsen calculus \cite{R} in the following way (keeping the same symbols for the knots):
a $1$-twist on $K_0$, a $1$-twist on $K_2$, and finally a $1$-twist on $K_1$.
This sequence of twisting gives the surgery description of Figure \ref{fig:surgery}.
\end{proof}

Let $\mu$ be the meridian of $K$ in $S^3$.

\begin{lemma}\label{lem:meridian}
If $\alpha\not\in \mathbb{Z}$ and $\beta\in\mathbb{Z}$ with respect to the original framing of $K$, 
then $\mu=1/0$ on $K_0$ in the surgery description of Figure \ref{fig:surgery}.
If $\alpha\in \mathbb{Z}$ and $\beta\not\in\mathbb{Z}$, then $\mu=1/1$.
\end{lemma}

\begin{proof}
Let $\mu=m/n$ on $K_0$ in Figure \ref{fig:surgery}.
By Lemma \ref{lem:surgery}, $\alpha=1/2$ and $\beta=3$ on $K_0$.
If $\alpha\not\in \mathbb{Z}$ and $\beta\in\mathbb{Z}$ with respect to the original framing of $K$, 
then $\Delta(\mu,\alpha)=|2m-n|=2$ and $\Delta(\mu,\beta)=|3n-m|=1$.
(Recall that $\alpha$ is half-integral.)
This gives $(m,n)=(1,0)$ or $(-1,0)$.  Hence $\mu=1/0$.

If $\alpha\in\mathbb{Z}$ and $\beta\not\in\mathbb{Z}$ with respect to the original framing of $K$, 
then a similar calculation shows that $(m,n)=(1,1)$ or $(-1,-1)$.
Thus $\mu=1/1$.
\end{proof}

Therefore, $K$ can be identified with $K_0$ in Figure \ref{fig:surgery} when $\alpha\not\in \mathbb{Z}$ and $\beta\in\mathbb{Z}$.
Otherwise, $K$ is the image of $K_0$ after a $(-1)$-twisting along $K_0$.

\begin{lemma}\label{lem:twist}
If $\alpha\not\in \mathbb{Z}$ and $\beta\in\mathbb{Z}$ with respect to the original framing of $K$, 
then $4p+3q=\pm 1$.
Otherwise, $4q-3p=\pm 1$.
\end{lemma}

\begin{figure}[tb]
%\blankbox{1.0\columnwidth}{2in}
\includegraphics*[scale=0.6]{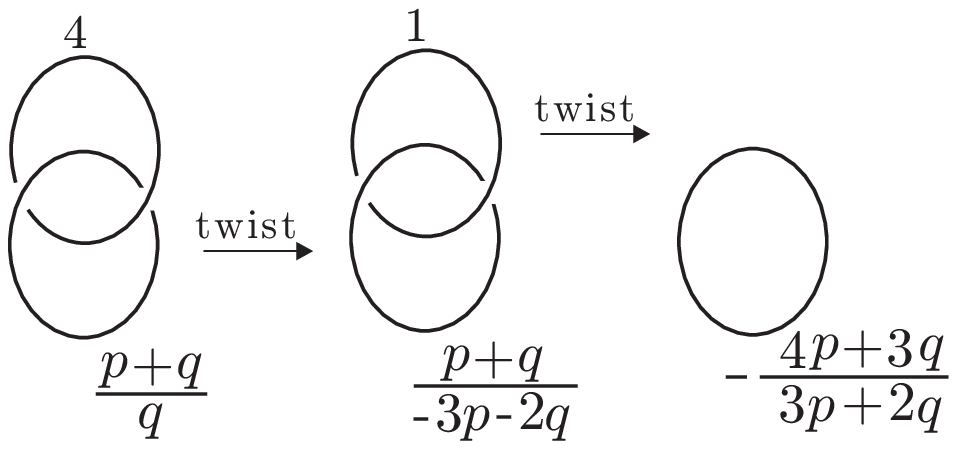}
\caption{}\label{fig:calculus}
\end{figure}

\begin{figure}[tb]
%\blankbox{1.0\columnwidth}{2in}
\includegraphics*[scale=0.6]{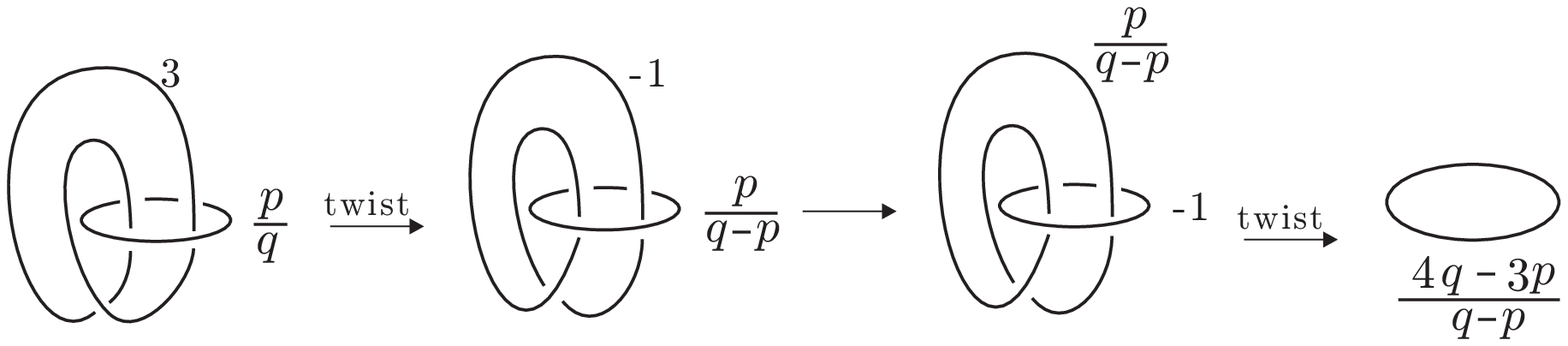}
\caption{}\label{fig:calculus2}
\end{figure}

\begin{proof}
Assume that $\alpha\not\in \mathbb{Z}$ and $\beta\in\mathbb{Z}$.
By Lemma \ref{lem:meridian}, the surgery description of Figure \ref{fig:surgery}, deleting $K_0$, gives $S^3$.
A Kirby-Rolfsen calculus shows that this is equivalent to the trivial knot with framing $-(4p+3q)/(3p+2q)$ as in Figure \ref{fig:calculus}.
Hence $4p+3q=\pm 1$.

When $\alpha\in \mathbb{Z}$ and $\beta\not\in\mathbb{Z}$,
first do a $(-1)$-twist on $K_0$ in Figure \ref{fig:surgery}, and eliminate $K_0$.
This gives $S^3$, and Figure \ref{fig:calculus2} shows that this is equivalent to the trivial knot with framing $(4q-3p)/(q-p)$.
Hence $4q-3p=\pm 1$.
\end{proof}

\begin{figure}[tb]
%\blankbox{1.0\columnwidth}{2in}
\includegraphics*[scale=0.56]{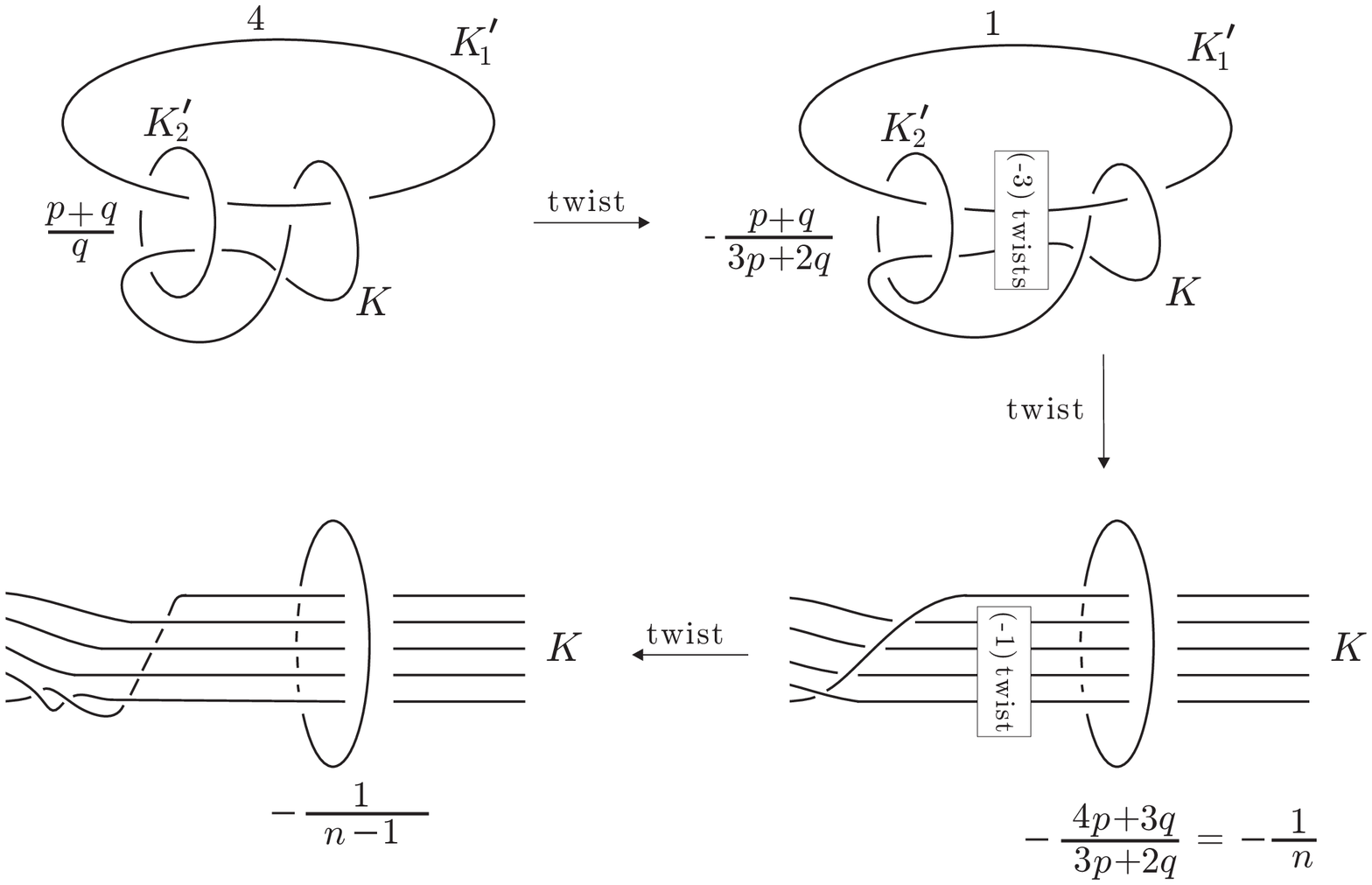}
\caption{}\label{fig:final}
\end{figure}

\begin{figure}[tb]
%\blankbox{1.0\columnwidth}{2in}
\includegraphics*[scale=0.56]{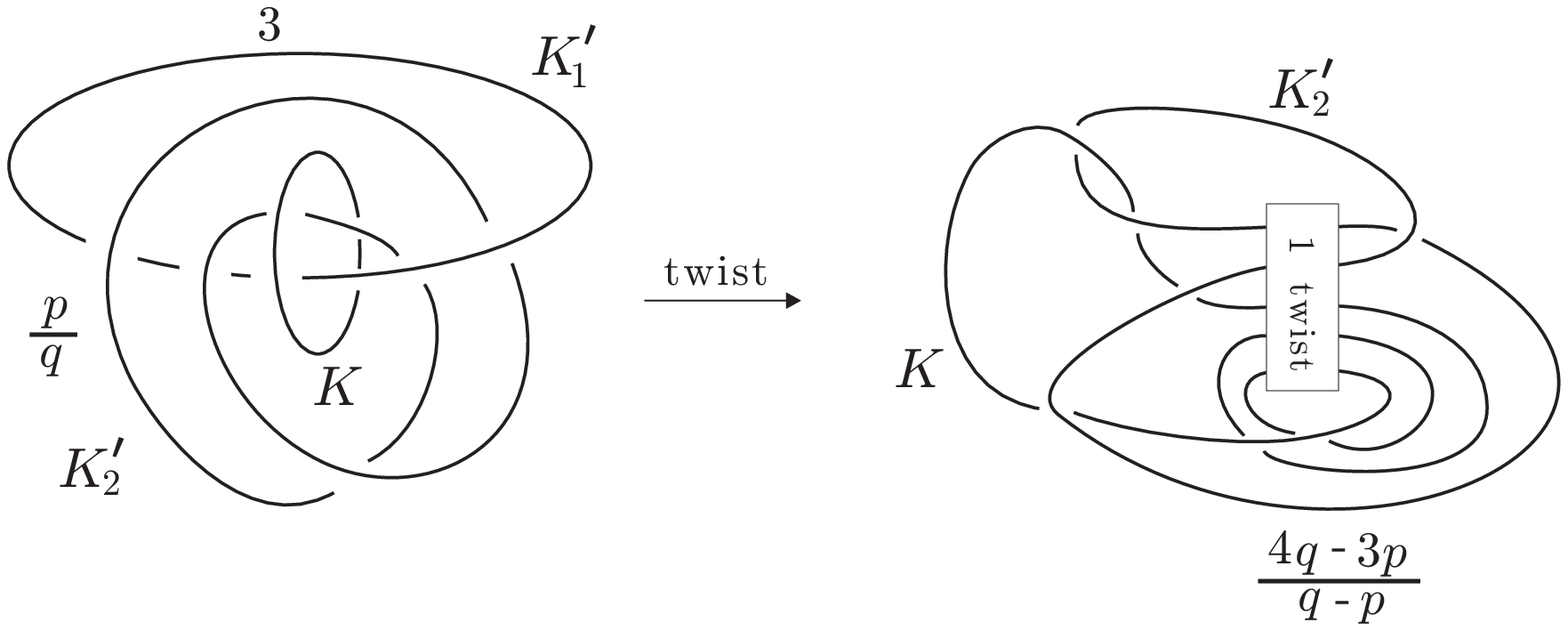}
\caption{}\label{fig:final2}
\end{figure}

\begin{proof}[Proof of Theorem \ref{thm:main}]
First assume that $\alpha\not\in \mathbb{Z}$ and $\beta\in\mathbb{Z}$.
We modify the surgery link of Figure \ref{fig:surgery} furthermore.
After a $(-3)$-twisting along $K_2'$, do a $(-1)$-twist along $K_1'$.
See Figure \ref{fig:final}, where $K=K_0$ is expressed in a braid form.
The framing on $K_2'$ is $-(4p+3q)/(3p+2q)=\pm 1/(3p+2q)$ by Lemma \ref{lem:twist}, and
$\alpha=-37/2$ and $\beta=-16$.

Put $n=(4p+3q)(3p+2q)=\pm (3p+2q)$.
After a $1$-twisting along $K_2'$, we have the Whitehead sister link, where
the framing on the unknotted component is $-1/(n-1)$.
Hence $K=k(2,-1,n,0)$, and $\alpha=25n-37/2$, $\beta=25n-16$.
It is easy to see that $n\ne 1$ and $n$ can be any other value.

Next, assume that $\alpha\in \mathbb{Z}$ and $\beta\not\in\mathbb{Z}$.
In this case, we perform a $(-1)$-twist along $K_0$ in Figure \ref{fig:surgery}.
Then $\mu=1/0$ in the resulting surgery link, and so $K$ can be identified with $K_0$.
Figure \ref{fig:final2} starts from this link.
At this point, $\alpha=1$ and $\beta=-3/2$.
After a $(-1)$-twist along $K_2'$, do a $1$-twist along $K_1'$.
Figure \ref{fig:final2} shows only the last result, where
the framing on $K_2'$ is $(4q-3p)/(q-p)=\pm 1/(q-p)$ by Lemma \ref{lem:twist},
and $\alpha=9$, $\beta=13/2$.
Again, this is the Whitehead sister link with the framing $\pm 1/(q-p)$ on the unknotted component.
Hence $K=k(2,-1,n,0)$, where $n=-(4q-3p)(q-p)+1=\pm (q-p)+1$, and $\alpha=25n-16$, $\beta=25n-37/2$.
\end{proof}

%%%%%%%%%%%%%%%%%%%%%%%%%%%%%%%%
\section*{Acknowledgements}
I would like to thank Professor Mario Eudave-Mu\~{n}oz for useful conversations.
I would also like to thank the referee for his helpful comments. 

%%%%%%%%%%%%%%%%%%%%%%%%%%%%%%%%%%
\bibliographystyle{amsplain}

\end{document}